\documentclass[11pt,reqno]{amsart}
\usepackage[T1]{fontenc}
\usepackage{lmodern}
\usepackage{amsmath,amssymb,mathtools,mathrsfs}
\usepackage{amsthm}
\usepackage{enumitem}
\usepackage{tikz-cd}
\usepackage[margin=1.12in]{geometry}
\usepackage{microtype}
\usepackage[colorlinks=true,linkcolor=blue,citecolor=blue,urlcolor=blue]{hyperref}
\hypersetup{pdftitle={Products and intersections of verbally prime T-ideals with applications to path ideals},pdfauthor={Yihao Zheng and Shenglin Zhu}}
\newtheorem{maintheorem}{Theorem}

\newcommand{\kk}{\Bbbk}
\newcommand{\KK}{\mathbb{K}}
\newcommand{\F}{\mathcal F}
\newcommand{\E}{\mathcal E}
\newcommand{\Id}{\operatorname{Id}}
\newcommand{\Mat}{\operatorname{M}}

\newcommand{\sgn}{\operatorname{sgn}}
\newcommand{\Ztwo}{\mathbb Z/2\mathbb Z}
\newcommand{\Qtype}{\mathsf Q}
\newcommand{\st}{\operatorname{St}}

\newtheorem{theorem}{Theorem}[section]
\newtheorem{lemma}[theorem]{Lemma}
\newtheorem{proposition}[theorem]{Proposition}
\newtheorem{corollary}[theorem]{Corollary}
\theoremstyle{definition}
\newtheorem{definition}[theorem]{Definition}
\newtheorem{example}[theorem]{Example}
\theoremstyle{remark}
\newtheorem{remark}[theorem]{Remark}
\numberwithin{equation}{section}
\setlist[enumerate]{label=\textup{(\roman*)},leftmargin=2.2em}
\setlist[itemize]{leftmargin=1.5em}
\title[Products and intersections of verbally prime $T$-ideals]{Inclusions and intersections of products
of verbally prime $T$-ideals}
\author{Yihao Zheng}
\address{School of Mathematical Sciences, Fudan University, Shanghai 200433, China}
\email{yhzheng24@m.fudan.edu.cn}
\author{Shenglin Zhu}
\address{School of Mathematical Sciences, Fudan University, Shanghai 200433, China}
\email{mazhusl@fudan.edu.cn}
\date{}
\begin{document}
\begin{abstract}
We study inclusions and intersections of products of verbally prime
$T$-ideals over a field of characteristic zero.  Our main tool is a Vandermonde-type construction of graded
polynomials with a uniformly bounded number of variables.
Together with a reduction to finite exterior algebras,
this construction yields separating polynomials that work
uniformly under bounds on Grassmann envelope representations.
We prove that, under these bounds, an intersection of
$T$-ideals is contained in such a product only if one member is.
This proves Berele's
Conjecture~3 and, consequently, his Conjecture~1 \cite{Berele2012}. We also characterize inclusions between ordered products
by an explicit combinatorial criterion, proving his Conjecture~2 and answering
his Question~2 \cite{Berele2012}.
These results characterize inclusions between arbitrary
intersections of products of nonzero proper verbally
prime $T$-ideals and show that every nonzero such
intersection has a unique finite irredundant representation.
As an application, we characterize finite quivers whose
path algebras have the same polynomial identities as some,
equivalently every, positive power of their arrow ideals.

\end{abstract}
\maketitle

% ===== parts/introduction =====
\section{Introduction}\label{sec:introduction}

Polynomial identities relate the structure of an associative
algebra to the relations satisfied by all its elements.
Their theory occupies an important place in noncommutative
algebra and has close connections with representation theory,
invariant theory, algebraic combinatorics, the theory of
division algebras, and noncommutative algebraic geometry;
see, for example,   
\cite{Bresar2014},  \cite{DrenskyFormanek2004},  \cite{Formanek1991},   
\cite{GiambrunoZaicev2005},  \cite{Jacobson1975},  \cite{ReichsteinVonessen2007}.

%Polynomial identities relate the structure of an associative
%algebra to the relations satisfied by all its elements. They
%provide invariants for comparing algebras and defining varieties,
%and their study connects the structure theory of algebras with
%representations of symmetric groups and asymptotic growth
%\cite{Bresar2014,GiambrunoZaicev2005}. The Amitsur--Levitzki
%theorem \cite{AmitsurLevitzki1950}, which determines the least
%degree of a nonzero identity of a full matrix algebra, is a
%basic example of how identities recover structural information.
%More generally, the inclusion of one identity ideal in another
%expresses a relation between the corresponding varieties of
%algebras. Determining such inclusions is therefore a natural
%part of understanding which structural information polynomial
%identities retain.

The polynomial identities of an algebra form a $T$-ideal,
so comparing the identities of different algebras leads
naturally to questions about inclusions between $T$-ideals.
Kemer's theory
\cite{Kemer1991} describes nonzero proper $T$-ideals as identity
ideals of Grassmann envelopes of finite-dimensional superalgebras
and classifies the verbally prime $T$-ideals. Products of these
ideals have a further role in the study of codimension growth.
Giambruno and Zaicev
\cite{GiambrunoZaicev2003},  \cite{GiambrunoZaicev2005}
study the structure and factorization of these products
in connection with minimal superalgebras and the codimension
growth of varieties of algebras.
Their work provides a basis for studying inclusions between
such products.

% Kemer's theory
%\cite{Kemer1991} describes nonzero proper $T$-ideals as identity
%ideals of Grassmann envelopes of finite-dimensional superalgebras
%and classifies the verbally prime $T$-ideals. Products of these
%ideals have a further role in the study of codimension growth.
%Giambruno and Zaicev
%\cite{GiambrunoZaicev2003,GiambrunoZaicev2005} relate them to
%minimal superalgebras and to varieties that are minimal for a
%given exponential rate of growth. They also established the
%uniqueness of the ordered factorization of a product of nonzero
%proper verbally prime $T$-ideals. These results give a structural
%reason to study products, as well as to ask how different
%products can be compared. 

Products and intersections also arise from concrete algebra
constructions. Lewin's theorem \cite{Lewin1974} accounts for
products of identity ideals in upper triangular matrix
constructions, whereas intersections describe the identities
common to the summands of a direct sum. In the setting of
incidence algebras, Berele \cite{Berele1984} expresses the
identity ideal as an intersection over maximal chains, with
each chain contributing a product of matrix identity ideals.
%Thus the combinatorics of a directed graph leads naturally to
%intersections of products of verbally prime $T$-ideals.

Berele \cite{Berele2012} gives partial answers to the question
of when two such intersections are equal or one is contained
in the other, but leaves the general problem open.
He formulates three conjectures concerning chain embeddings,
inclusions between verbally prime $T$-ideals, and finite
intersections.
In this paper, we prove stronger results on intersections of
$T$-ideals, inclusions between products of verbally prime
$T$-ideals, and inclusions between arbitrary intersections
of such products.
 These results imply all three
conjectures and answer his Question~2.

Our first result is an intersection theorem for possibly
infinite families of $T$-ideals whose Grassmann envelope
representations satisfy uniform bounds.

%He proved that a product of proper verbally
%prime $T$-ideals cannot be expressed as a finite intersection
%of strictly larger $T$-ideals. His Conjecture~3 asks for the
%stronger conclusion that containment of a finite intersection
%in such a product already forces containment of one member.
%This is precisely the additional property needed to pass from
%comparisons of individual products to comparisons of their
%intersections.

% Throughout,
%$\kk$ is a field of characteristic zero, and $T$-ideals are
%considered in the free associative algebra
%$\F=\kk\langle x_1,x_2,\cdots\rangle_+$ without constant term.
%We use $\E(B)$ for the Grassmann envelope of a superalgebra
%$B$ and $J(B)$ for its Jacobson radical.

\begin{maintheorem}[Theorem~\ref{thm:bounded-intersection}]
\label{main:intersection}
Let $\{T_\lambda\}_{\lambda\in\Lambda}$ be a nonempty
family of $T$-ideals. Suppose that there are integers
$s_0\geq0$ and $\nu_0\geq1$ such that every nonzero
proper $T_\lambda$ has the form
$T_\lambda=\Id_{\kk}(\E_{\KK_\lambda}(B_\lambda))$,
where $\KK_\lambda$ is an algebraically closed extension
of $\kk$ and $B_\lambda$ is a finite-dimensional
$\KK_\lambda$-superalgebra satisfying
$\dim_{\KK_\lambda}(B_\lambda/J(B_\lambda))\leq s_0$
and $J(B_\lambda)^{\nu_0}=0$.
If $P_1,\cdots,P_r$ are proper verbally prime
$T$-ideals, where $r\geq1$, then
$\bigcap_{\lambda\in\Lambda}T_\lambda
\subseteq P_1\cdots P_r$ only if
$T_\lambda\subseteq P_1\cdots P_r$
for some $\lambda\in\Lambda$.
\end{maintheorem}

For a finite family, Kemer's representability theorem \cite{Kemer1987}  supplies
the required finite-dimensional superalgebras, and finiteness
of the family gives common bounds. Theorem~\ref{main:intersection}
therefore proves Conjecture~3
(Theorem~\ref{thm:finite-intersection}).
Berele \cite[Corollary~2.4]{Berele2012} shows that
Conjecture~3 implies Conjecture~1.
Thus we also obtain Conjecture~1
(Corollary~\ref{cor:berele-one}).

To apply Theorem~\ref{main:intersection} to intersections of
products, one must next decide when one product is contained
in another.
Berele \cite{Berele2012} characterizes inclusions
between products of ordinary matrix identity ideals
combinatorially. For finite-dimensional graded-simple
algebras graded
by an abelian group over an algebraically closed
field of characteristic zero, David \cite{David2012}
characterizes inclusions of graded identity ideals by
graded embeddings in the opposite direction.
Our second result gives a criterion  for  inclusion between  products
of ordinary identity ideals of Grassmann envelopes, allowing
both types of graded-simple superalgebras.

%For a nonempty ordered family
%$\mathcal S=(S_1,\cdots,S_r)$ of split graded-simple
%superalgebras, let
%$I(\mathcal S)=P(S_1)\cdots P(S_r)$, where
%$P(S)=\Id(\E(S))$.
%The relation $\preceq_{\mathrm{sup}}$ in
%Definition~\ref{def:super-capacity} is specified by a
%nondecreasing map between the index sets and explicit
%inequalities on the matrix and parity parameters.

\begin{maintheorem}[Theorem~\ref{thm:mixed-capacity}]
\label{main:products}
Let $\mathcal S$ and $\mathcal R$ be nonempty ordered families of split graded-simple superalgebras over a field of characteristic zero. Then
$
I(\mathcal R)\subseteq I(\mathcal S)$ if and only if
 $\mathcal S\preceq_{\mathrm{sup}}\mathcal R$.
\end{maintheorem}

Taking both ordered families to have one member proves
Berele's Conjecture~2
(Corollary~\ref{cor:prime-inclusions}).
Taking $\mathcal R$ to have one member answers his
Question~2 (Corollary~\ref{cor:question-two}), which asks when a verbally prime $T$-ideal
is contained in a product of such ideals.

The preceding results also allow us to compare intersections
and determine which products are needed to represent them.
This extends the uniqueness question from the factors of
a single product to the products occurring in an intersection.
The index set of the intersection need not be finite.
%For the next result, write every matrix-type member as
%$\Mat_{a|b}(\kk)$ with $a\geq b\geq0$ and $a+b>0$,
%and every other member as $\Qtype_n$ with $n\geq1$.
%For a nonempty set $\Gamma$ of nonempty finite ordered
%families of these superalgebras, let
%$\mathcal I_\Gamma=\bigcap_{\gamma\in\Gamma}I(\gamma)$.
%When $\Gamma$ is finite, let $\Gamma_{\max}$ denote its
%maximal members with respect to $\preceq_{\mathrm{sup}}$.

\begin{maintheorem}[Theorem~\ref{thm:intersection-capacity} and Theorem~\ref{thm:unique-antichain}]
\label{main:intersections}
Let $\Gamma$ and $\Delta$ be nonempty sets of nonempty
finite ordered families of split graded-simple
superalgebras in the fixed normalized forms.
\begin{enumerate}
\item
$\mathcal I_\Gamma\subseteq\mathcal I_\Delta$
if and only if $\Gamma$ is infinite or, for every
$\delta\in\Delta$, there is $\gamma\in\Gamma$ with
$\delta\preceq_{\mathrm{sup}}\gamma$.
\item
If $\Gamma$ is finite, then
$\mathcal I_\Gamma=\mathcal I_{\Gamma_{\max}}$.
Moreover, $\mathcal I_\Gamma=\mathcal I_\Delta$ if and
only if either both sets are infinite, or both are
finite and $\Gamma_{\max}=\Delta_{\max}$.
\end{enumerate}

In particular, let $\Lambda$ be a nonempty set, and for
each $\lambda\in\Lambda$, let
$T_\lambda=P_{\lambda,1}\cdots P_{\lambda,r_\lambda}$,
where $r_\lambda\geq1$ and each $P_{\lambda,i}$ is a
nonzero proper verbally prime $T$-ideal.
Let $\mathfrak a=\bigcap_{\lambda\in\Lambda}T_\lambda$.

If $\mathfrak a\ne0$, then there is a finite nonempty
subset $\Lambda_0\subseteq\Lambda$ such that
$\mathfrak a=\bigcap_{\lambda\in\Lambda_0}T_\lambda$
and
$\mathfrak a\subsetneq
\bigcap_{\mu\in\Lambda_0\setminus\{\lambda\}}T_\mu$
for every $\lambda\in\Lambda_0$.
The set $\{T_\lambda\mid\lambda\in\Lambda_0\}$
is uniquely determined by $\mathfrak a$.

If $\mathfrak a=0$, then
$\bigcap_{\lambda\in\Lambda\setminus\{\lambda_0\}}
T_\lambda=0$
for every $\lambda_0\in\Lambda$.
Here the empty intersection is defined to be $\F$.

\end{maintheorem}

The comparison of identity ideals also applies to path
algebras, whose polynomial identities are studied by
Cerulli Irelli, De Loera Ch\'avez, and Pascucci
\cite{CerulliIrelliEtAl2026}.
They characterize finite quivers whose path algebras are
PI and determine the identity ideals in the acyclic case.
Berele et al.~\cite{BereleEtAl2026} show that every
subalgebra generated by paths in a PI quiver has the same
identities as an associated finite-dimensional incidence
algebra.
These results bring the comparison of path algebras
and path ideals into the setting of products and
intersections considered here.

The following result gives a combinatorial criterion for
equality between the identity ideals of $\kk Q$ and a
positive power of its arrow ideal, and shows that equality
for one positive power implies equality for all positive
powers.

%For a finite quiver $Q$, let $\kk Q_{\ge m}$ be the span
%of paths of length at least $m$. These ideals are the
%positive powers of the arrow ideal $\kk Q_{\ge1}$.
%We ask when passing to these ideals leaves the identity
%ideal of $\kk Q$ unchanged.
 % Here $\mathcal P(Q)$ is the partially
%ordered set of strongly connected components (SCCs)
%under reachability, $\mathbf w(\mathcal C)$ records
%the numbers of vertices in the SCCs of a chain, and
%$\preceq$ is the relation in Definition~\ref{def:chain-order}.

\begin{maintheorem}[Corollary~\ref{cor:all-tail-equality}]
\label{main:paths}
Let $Q$ be a finite quiver. If $Q$ is PI, the following
conditions are equivalent. 
\begin{enumerate}
\item $\Id(\kk Q)=\Id(\kk Q_{\ge m})$ for some $m\ge1$.
\item $\Id(\kk Q)=\Id(\kk Q_{\ge m})$ for every $m\ge1$.
\item For every maximal chain $\mathcal C$ in $\mathcal P(Q)$,
there is a chain $\mathcal D$ which is maximal among chains
consisting entirely of cyclic SCCs and satisfies
$\mathbf w(\mathcal C)\preceq\mathbf w(\mathcal D)$.
\end{enumerate}
If $Q$ is not PI, $\Id(\kk Q)=\Id(\kk Q_{\ge m})=0$ for every $m\ge1$.
\end{maintheorem}

The main tool in the proofs is a Vandermonde-type construction
of graded polynomials associated with an ordered family
$\mathcal S$ of graded-simple superalgebras.
Their degrees can be increased while the number of variables
remains bounded solely in terms of $\mathcal S$.
This differs from the construction of the strongly full
polynomials used by Aljadeff and Karasik
\cite{AljadeffKarasik2024}, which relies on alternation and
symmetrization.
For a finite-dimensional superalgebra over an algebraically
closed field, nonzero graded evaluations of these polynomials
of sufficiently large degree yield a minimal superalgebra
whose graded-simple subalgebras are graded-isomorphic to
the members of $\mathcal S$ in the given order.

The bound on the number of variables allows us to reduce
to graded tensor products with a finite exterior algebra
whose number of generators is independent of the polynomial
degree. We can therefore choose the degree sufficiently
large relative to the nilpotency index of the radical
without enlarging the exterior algebra.
Together with the product formula for minimal superalgebras,
this gives a single multilinear polynomial that determines
containment in $I(\mathcal S)$ uniformly under common bounds
on the dimensions of the semisimple quotients and the
nilpotency indices of the radicals.
This uniform choice yields Theorem~\ref{main:intersection}.
The minimal subalgebras obtained from nonzero evaluations
also yield the index requirements in
Theorem~\ref{main:products}.

%The proofs of the general results use separating polynomials
%associated with ordered families of graded-simple
%superalgebras. For sufficiently large values of the parameter,
%a nonzero graded evaluation of any of these polynomials  yields a
%minimal superalgebra with the given graded-simple
%subalgebras. A reduction to
%graded tensor products with finite exterior algebras
%makes this construction uniform under the bounds in
%Theorem~\ref{main:intersection}. Together with the product
%formula for minimal superalgebras, this gives both the
%intersection theorem and the inclusion criterion.

The paper is organized as follows.
Section~\ref{sec:preliminaries} collects the preliminaries.
Section~\ref{sec:grassmann} establishes the reductions to graded identities
and finite exterior algebras.  Section~\ref{sec:separation} constructs the separating
polynomials and the minimal superalgebras arising from
their nonzero evaluations.
Section~\ref{sec:finite-intersections} proves the intersection
theorems and deduces Berele's Conjecture~3 and Conjecture~1.
Section~\ref{sec:capacity} gives the criterion for inclusion
between products, proves Berele's Conjecture~2, and answers Question~2.
Section~\ref{sec:antichains} characterizes inclusions between
intersections of products and proves that every nonzero
intersection has a unique finite irredundant representation.
Section~\ref{sec:path-applications} applies these results
to path algebras and path ideals.
% and also discusses
%stabilization, numerical criteria, and examples.

\medskip
\noindent\textbf{Acknowledgments.}
The authors thank Giovanni Cerulli Irelli for helpful discussions and  comments.

\medskip
\noindent\textbf{AI Declaration.}
The proofs of Theorem~\ref{thm:bounded-intersection} and
Theorem~\ref{thm:mixed-capacity} were found by ChatGPT 6.
The authors' contribution to these proofs is to verify, organize,
and present them in a rigorous form for a human audience.
ChatGPT 6 was additionally used for language refinement,
proofreading, and checking the presentation of some parts of
the manuscript.
The authors have independently verified the mathematical content
and take full responsibility for the exposition and correctness
of the manuscript.

% ===== parts/preliminaries =====
\section{Preliminaries}\label{sec:preliminaries}

Throughout the paper, $\kk$ is a field of characteristic zero. Algebras are associative and need not  be unital.  If  \(A\) is a unital algebra, a subalgebra $C$ of $A$ need not be unital, and even if \(C\) is unital,   \(1_C\) need not equal \(1_A\). Subalgebras generated by subsets of \(A\) are taken without adjoining \(1_A\). An algebra embedding is an injective algebra homomorphism, with no requirement that identity elements exist or be preserved. 
We recall the basic definitions and results concerning
polynomial identities and superalgebras. We also recall the basic concepts concerning quivers and 
path algebras. Standard references are  
\cite{AssemSimsonSkowronski2006},  \cite{BangJensenGutin2009},  \cite{Bresar2014},  
\cite{DerksenWeyman2017},  \cite{GiambrunoZaicev2005},  \cite{Kemer1991}.

\subsection{Polynomial identities}\label{subsec:polynomial-identities}

We fix the conventions for polynomial identities and verbally prime
$T$-ideals, and recall the degree bound for identities of matrix algebras.
Let $\F=\kk\langle x_1,x_2,\cdots\rangle_+$ be the free
associative algebra without constant term.

\begin{definition}\label{def:polynomial-identities}
A polynomial $f\in\F$ is an identity of an algebra $A$ if
$f$ vanishes under every evaluation in $A$.
Let $\Id(A)$ denote the ideal of these identities.
The algebra $A$ is PI if $\Id(A)\ne0$.
A $T$-ideal is a two-sided ideal of $\F$ invariant under every
endomorphism of $\F$. In particular, $\Id(A)$ is a $T$-ideal
and $\Id(0)=\F$.
\end{definition}

Throughout, for any field $L$ of characteristic zero,
all $T$-ideals, including the identity ideals of unital
algebras, are considered in
$L\langle x_1,x_2,\cdots\rangle_+$.
Indeed, every proper $T$-ideal of
$L\langle x_1,x_2,\cdots\rangle$
is contained in $L\langle x_1,x_2,\cdots\rangle_+$
and is a $T$-ideal there.
An evaluation in a unital algebra $B$ may assign $1_B$
to any variable.

Let $\KK/\kk$ be a field extension and $A$ a $\KK$-algebra.
We denote by $\Id_{\KK}(A)$ the $T$-ideal of polynomial
identities of $A$ in
$\KK\langle x_1,x_2,\cdots\rangle_+$,
and define $
\Id_{\kk}(A)=\Id_{\KK}(A)\cap\F$.

\begin{definition}\label{def:verbally-prime}
A proper $T$-ideal $P$ is verbally prime if, for any $T$-ideals
$U$ and $V$, the inclusion $UV\subseteq P$ implies
$U\subseteq P$ or $V\subseteq P$.
\end{definition}

\begin{definition}\label{def:standard-polynomial}
For $d\ge1$, the standard polynomial
$\st_d(x_1,\cdots,x_d)$ of degree $d$ is defined by
$\st_d=\sum_{\sigma\in\mathfrak S_d}
\sgn(\sigma)x_{\sigma(1)}\cdots x_{\sigma(d)}$,
where $\mathfrak S_d$ is the symmetric group on
$\{1,\cdots,d\}$.
\end{definition}

For $n\ge1$, let $I_n=\Id(\Mat_n(\kk))$. Amitsur and Levitzki \cite{AmitsurLevitzki1950}  prove the   celebrated  Amitsur--Levitzki theorem, which is stated as follows.
\begin{theorem}\label{thm:amitsur-levitzki}
For every $n\ge1$, we have $\st_{2n}\in I_n$, and every nonzero
polynomial in $I_n$ has  degree at least $2n$.
\end{theorem}

\subsection{Superalgebras and Grassmann envelopes}\label{subsec:superalgebras-envelopes}

We recall graded identities and complete linearization, and then
define Grassmann envelopes and graded tensor products with exterior
algebras.

\begin{definition}\label{def:superalgebras}
A superalgebra is an algebra $B=B_0\oplus B_1$ graded by
$\Ztwo$, so that $B_iB_j\subseteq B_{i+j}$ for $i,j\in\Ztwo$.
For a nonzero homogeneous element $b\in B_i$, its parity is
$|b|=i$.
The superalgebra $B$ is supercommutative if
$
bc=(-1)^{|b||c|}cb
$
for all homogeneous $b,c\in B$.
A nonzero superalgebra $S$ is graded-simple if $S^2\ne0$ and
its only graded two-sided ideals are $0$ and $S$.
\end{definition}

All gradings in this paper are $\Ztwo$-gradings unless
explicitly stated otherwise.

\begin{definition}\label{def:graded-identities}
Let $B$ be a superalgebra, and let $f(x_1,\cdots,x_n)$ have
assigned variable parities
$\delta=(\delta_1,\cdots,\delta_n)\in(\Ztwo)^n$.
The polynomial $f$ is a graded identity of $B$ if
$
f(b_1,\cdots,b_n)=0$ for all $b_i\in B_{\delta_i}$.
A graded $T$-ideal is a two-sided ideal of the free graded
algebra invariant under every degree-preserving endomorphism.
\end{definition}

We introduce complete linearization, which allows us to work
with multilinear identities in both the ordinary and graded settings.

\begin{definition}\label{def:linearization}
Let $f(x_1,\cdots,x_n)$ be multihomogeneous of degree $d_i\geq1$
in each variable $x_i$. Introduce distinct variables
$x_{i1},\cdots,x_{i,d_i}$ for each $i$.
The complete linearization $\widetilde f$ of $f$ is the
component of
$
f(\sum_{j=1}^{d_1}x_{1j},\cdots,
\sum_{j=1}^{d_n}x_{nj})
$
that has degree one in every variable $x_{ij}$.
For a graded polynomial, each $x_{ij}$ has the parity of $x_i$.
We do not divide by factorials in this definition.
\end{definition}

The following standard result follows from complete linearization;
see \cite{GiambrunoZaicev2005}
for the ordinary case. The graded case follows by the same argument.

\begin{theorem}\label{thm:linearization}
Every $T$-ideal is generated, as a $T$-ideal, by its
multilinear polynomials.
A multihomogeneous polynomial $f$ and its complete
linearization $\widetilde f$ generate the same $T$-ideal.
The same statements hold for graded $T$-ideals, provided
each new variable has the same parity as the variable
it replaces.
\end{theorem}

Let $\E$ be the exterior algebra on countably many generators,
graded by the parity of exterior degree. Let $\E^{(N)}$ be the
exterior algebra on $N$ generators, where $\E^{(0)}=\kk$.
We write $\E_{\KK}$ and $\E_{\KK}^{(N)}$ to indicate that these
exterior algebras are defined over $\KK$.
The even part $\E_0$ is spanned by $1$ and products of an even number of
distinct exterior generators, and the odd part $\E_1$ by products of an odd
number. These exterior algebras are supercommutative in the
sense of Definition~\ref{def:superalgebras}.

\begin{definition}\label{def:grassmann-envelope}
The Grassmann envelope of a superalgebra $B$ is
\[
\E(B)=(B_0\otimes\E_0)\oplus(B_1\otimes\E_1),
\]
with multiplication inherited from the ordinary tensor product.
The graded tensor product $B\mathbin{\widehat{\otimes}}\E^{(N)}$ has
$|b\otimes\eta|=|b|+|\eta|$ and multiplication
\[
(b\otimes\eta)(c\otimes\theta)
=(-1)^{|\eta||c|}bc\otimes\eta\theta
\]
for homogeneous elements. The same conventions define
$B\mathbin{\widehat{\otimes}}\E$. 
\end{definition}

\subsection{Structure of superalgebras and representability of
\texorpdfstring{$T$}{T}-ideals}
\label{subsec:graded-simple-structure}

We recall structural results on finite-dimensional superalgebras
and the representation of nonzero proper $T$-ideals as identity
ideals of Grassmann envelopes.

\begin{definition}\label{def:split-simple}
For $p,q\geq0$ with $p+q>0$, let $\Mat_{p|q}(\kk)$ be the
matrix algebra $\Mat_{p+q}(\kk)$ with grading
\[
\begin{aligned}
\Mat_{p|q}(\kk)_0
&=\left\{\begin{pmatrix}A&0\\0&D\end{pmatrix}
\,\middle|\,A\in\Mat_p(\kk),\ D\in\Mat_q(\kk)\right\},\\
\Mat_{p|q}(\kk)_1
&=\left\{\begin{pmatrix}0&B\\C&0\end{pmatrix}
\,\middle|\,B\in\Mat_{p\times q}(\kk),\
C\in\Mat_{q\times p}(\kk)\right\}.
\end{aligned}
\]
If $p=0$ or $q=0$, the grading is concentrated in degree zero.

For $n\geq1$, let
\[
\Qtype_n=\Mat_n(\kk)\oplus u\Mat_n(\kk),
\]
where $u$ is odd, $u^2=1$, and $u$ commutes with $\Mat_n(\kk)$.
Its even part is $\Mat_n(\kk)$ and its odd part is $u\Mat_n(\kk)$.
A finite-dimensional graded-simple superalgebra over $\kk$ is
split if it is graded-isomorphic to one of these algebras.
\end{definition}

We recall the classification of finite-dimensional graded-simple
superalgebras over an algebraically closed field of characteristic zero; see
\cite{AljadeffKanelBelov2010} and
\cite{GiambrunoZaicev2005}.

\begin{theorem}\label{thm:graded-simple-classification}
Over an algebraically closed field of characteristic zero,
every finite-dimensional graded-simple superalgebra is
graded-isomorphic to one of the superalgebras in
Definition~\ref{def:split-simple}.
\end{theorem}

The following theorem gives Kemer's classification of verbally
prime $T$-ideals; see \cite{GiambrunoZaicev2005} and \cite{Kemer1991}.

\begin{theorem}\label{thm:kemer-classification}
The nonzero proper verbally prime $T$-ideals of $\F$ are
precisely the ideals $\Id(\E(S))$, where $S$ is a split
graded-simple superalgebra over $\kk$.
\end{theorem}

For arbitrary nonzero proper $T$-ideals, Kemer's representability
theorem \cite{AljadeffKanelBelov2010} gives the following result.

\begin{theorem}\label{thm:kemer-representability}
For every nonzero proper $T$-ideal $T\subseteq\F$, there are
an algebraically closed extension $\KK/\kk$ and a finite-dimensional
$\KK$-superalgebra $B$ such that
$
T=\Id_{\kk}(\E_{\KK}(B))$.
\end{theorem}

Let $J(B)$ denote the Jacobson radical of a finite-dimensional
superalgebra $B$. We also use the graded Wedderburn--Malcev
decomposition; see \cite{AljadeffKanelBelov2010}.

\begin{theorem}\label{thm:graded-wedderburn}
Let $B$ be a finite-dimensional superalgebra over an algebraically
closed field of characteristic zero. The radical $J(B)$ is graded,
and there is a graded semisimple subalgebra $S$ such that $
B=S\oplus J(B)$
as a vector space. 
\end{theorem}

\subsection{Minimal superalgebras and scalar extension}\label{subsec:minimal-scalar-extension}

We recall the definition of minimal superalgebras and the product
formula for the identity ideals of their Grassmann envelopes.
We also record how these identity ideals behave under scalar extension.
We first fix notation for this product.

\begin{definition}\label{def:ordered-family}
For a split graded-simple superalgebra $S$, let
$P(S)=\Id(\E(S))$.
For a nonempty ordered family $\mathcal S=(S_1,\cdots,S_r)$
of split graded-simple superalgebras, let
\[
I(\mathcal S)=P(S_1)\cdots P(S_r).
\]
Repeated members are allowed.
\end{definition}

\begin{definition}\label{def:minimal-superalgebra}
Let $C$ be a finite-dimensional superalgebra over an algebraically
closed field, with a graded Wedderburn--Malcev decomposition $
C=(S_1\oplus\cdots\oplus S_r)\oplus J(C)$,
where the $S_i$ are graded-simple, and let
$\mathcal S=(S_1,\cdots,S_r)$.
A nonzero homogeneous idempotent is minimal if it is not a sum
of two nonzero orthogonal homogeneous idempotents.
We say that $C$ is minimal with respect to $\mathcal S$ if there
are minimal nonzero homogeneous idempotents $e_i\in S_i$
and homogeneous elements
$
w_i\in e_iJ(C)e_{i+1}$ for $1\leq i<r$
such that $w_1\cdots w_{r-1}\ne0$ and $J(C)$ is generated by
the $w_i$ as an ideal of $C$. For $r=1$, this means $C=S_1$.
\end{definition}

A minimal homogeneous idempotent need not be primitive in the
underlying ungraded algebra. Under the graded embedding
\[
\Qtype_n\hookrightarrow\Mat_{n|n}(\kk),\qquad
X+uY\longmapsto\begin{pmatrix}X&Y\\Y&X\end{pmatrix},
\]
the standard minimal homogeneous diagonal idempotents are
$E_{ii}+E_{n+i,n+i}$ for $1\leq i\leq n$, where $E_{ab}$ denotes
a matrix unit of $\Mat_{2n}(\kk)$.

The following theorem expresses the identity ideal of the
Grassmann envelope of a minimal superalgebra as a product
of verbally prime $T$-ideals \cite{GiambrunoZaicev2003}.

\begin{theorem}\label{thm:minimal-product}
Let $C$ be minimal with respect to the ordered family
$\mathcal S=(S_1,\cdots,S_r)$ in
Definition~\ref{def:minimal-superalgebra}. Then $\Id(\E(C))=I(\mathcal S)$.
\end{theorem}

The algebras in Theorem~\ref{thm:kemer-representability} may be
defined over different algebraically closed extensions of $\kk$.
The next lemma gives formulas for scalar extension of identity
ideals of Grassmann envelopes and of products of such ideals.
 For a $\kk$-subspace $U\subseteq\F$, let ${\KK}U$ denote its
$\KK$-linear span in $\KK\langle x_1,x_2,\cdots\rangle_+$.

\begin{lemma}\label{lem:scalar-extension}
For every finite-dimensional $\kk$-superalgebra $A$ and every
extension $\KK/\kk$,
\[
\Id_{\KK}(\E_{\KK}(A\otimes_\kk \KK))=\KK\Id_\kk(\E_\kk(A)).
\]
In particular, let $\mathcal S=(S_1,\cdots,S_r)$ be an ordered
family of split graded-simple superalgebras over $\kk$, and let
$P_i=\Id_\kk(\E_\kk(S_i))$. Then
\[
\Id_{\KK}(\E_{\KK}(S_i\otimes_\kk \KK))={\KK}P_i,
\qquad
({\KK}P_1)\cdots({\KK}P_r)=\KK(P_1\cdots P_r).
\]
Consequently,
\[
\bigl(({\KK}P_1)\cdots({\KK}P_r)\bigr)\cap\F=P_1\cdots P_r.
\]
\end{lemma}

\begin{proof}
There is a natural isomorphism $
\E_{\KK}(A\otimes_\kk \KK)\cong\E_\kk(A)\otimes_\kk \KK$.
By \cite[Theorem~6.29]{Bresar2014}, every polynomial identity
of $\E_\kk(A)$ remains an identity after scalar extension.
Conversely, write a polynomial identity over $\KK$ as
$f=\sum_j\lambda_j f_j$, where the $\lambda_j\in \KK$ are
linearly independent over $\kk$ and the $f_j$ have
coefficients in $\kk$. Evaluating on $\E_\kk(A)$ shows
that each $f_j$ is an identity of $\E_\kk(A)$.
This proves
$
\Id_{\KK}(\E_{\KK}(A\otimes_\kk \KK))
=\KK\Id_\kk(\E_\kk(A))$.
The remaining assertions are trivial.
\end{proof}

\subsection{Quivers and path algebras}\label{subsec:quiver}

We  recall  some concepts concerning quivers and 
path algebras.
Throughout this paper, quivers are finite and their vertex sets
are nonempty. Loops, parallel arrows, and isolated vertices are allowed.

\begin{definition}\label{def:quiver-paths}
A quiver is a quadruple $Q=(Q_0,Q_1,s,t)$, where $Q_0$ is the
vertex set, $Q_1$ is the arrow set, and $s,t:Q_1\to Q_0$ assign
to each arrow its starting and terminal vertices.
A path of length $d\ge1$ is a sequence $a_1\cdots a_d$ of arrows
such that $t(a_i)=s(a_{i+1})$ for $1\le i<d$.
Its starting and terminal vertices are $s(a_1)$ and $t(a_d)$.
For each vertex $v$, let $e_v$ be the trivial path  at $v$.
A closed path has positive length and equal starting and terminal
vertices. A directed cycle on $n$ vertices has arrows
$1\to2\to\cdots\to n\to1$ and no other arrows; the case $n=1$
is a loop.

The path algebra $\kk Q$ has all paths as a basis.
The product $pq$ is their concatenation if the terminal vertex
of $p$ is the starting vertex of $q$, and is zero otherwise.
Its identity element is $\sum_{v\in Q_0}e_v$.
For $m\ge1$, let $\kk Q_{\ge m}$ be the span of paths of length
at least $m$. This is a two-sided ideal, regarded as a nonunital
algebra. The ideal $\kk Q_{\ge1}$ is the arrow ideal, and
$\kk Q_{\ge m}=(\kk Q_{\ge1})^m$.
A quiver is PI if its path algebra is PI.
\end{definition}

\begin{definition}\label{def:quiver-components}
Two vertices belong to the same strongly connected component (or SCC, for short)
if each is reachable from the other by a path.
Let $\mathcal P(Q)$ be the set of SCCs of $Q$.
For $C,D\in\mathcal P(Q)$, write $C\leq_Q D$ if a path, possibly
trivial, joins a vertex of $C$ to a vertex of $D$.
This is a partial order. Write $C<_Q D$ if $C\leq_Q D$ and $C\ne D$.
Let $w(C)$ be the number of vertices of $C$.
An SCC is cyclic if it contains a closed path.
A cyclic chain of SCCs is a nonempty chain with respect to $\leq_Q$ consisting entirely of cyclic SCCs.
\end{definition}

% ===== parts/grassmann =====
\section{Grassmann reductions for polynomial identities}\label{sec:grassmann}

We relate ordinary identities of Grassmann envelopes to graded
identities of graded tensor products with exterior algebras.
We then bound the number of exterior generators needed for
polynomials in a fixed number of variables, independently of
their degrees. We also bound the nilpotency index of the radical
after taking a graded tensor product with a finite exterior algebra.
These bounds are used to obtain uniform criteria for containment
of identity ideals of Grassmann envelopes in products of verbally
prime $T$-ideals.

\subsection{Ordinary and graded identities}\label{subsec:ordinary-graded-identities}

We use a sign transformation to relate ordinary identities of
Grassmann envelopes to graded identities of graded  tensor products with
exterior algebras. We then deduce that taking a graded tensor product
with a finite exterior algebra preserves the identities of the
Grassmann envelope.

\begin{definition}\label{def:sign-transformation}
Let $
f(x_1,\cdots,x_n)
=\sum_{\sigma\in\mathfrak S_n}
c_\sigma x_{\sigma(1)}\cdots x_{\sigma(n)}$
be a multilinear polynomial. 
For $\sigma\in\mathfrak S_n$, let
$
\mathcal G_\sigma
=
\{(i,j)\mid 1\leq i<j\leq n,\;
\sigma^{-1}(i)>\sigma^{-1}(j)\}$.
For a parity vector
$\delta=(\delta_1,\cdots,\delta_n)\in(\Ztwo)^n$, define
\[
\tau_\delta f
=
\sum_{\sigma\in\mathfrak S_n}
\chi_\delta(\sigma)c_\sigma
x_{\sigma(1)}\cdots x_{\sigma(n)},
\qquad\text{where}\qquad
\chi_\delta(\sigma)
=
(-1)^{\displaystyle
\sum\limits_{(i,j)\in\mathcal G_\sigma}
\delta_i\delta_j}.
\]
The exponent defining $\chi_\delta(\sigma)$ counts, modulo two, the pairs of odd variables
whose order in the product is opposite to their order in
$x_1\cdots x_n$. Thus $\tau_\delta^2=\mathrm{id}$.
\end{definition}

\begin{lemma}\label{lem:sign-transfer}
Let $B$ be a superalgebra and let $g$ be multilinear with variable parities $\delta$. Then $g$ is a graded identity of $B\mathbin{\widehat{\otimes}}\E$ if and only if $\tau_\delta g$ is an ordinary identity of $\E(B)$.
\end{lemma}

\begin{proof}
Let $h=\tau_\delta g$.
By \cite[Lemma~3.7.4]{GiambrunoZaicev2005} and multilinearity,
$h\in\Id(\E(B))$ if and only if
$
(\tau_\epsilon h)(b_1,\cdots,b_n)=0$
for every $\epsilon\in(\Ztwo)^n$ and all
$b_i\in B_{\epsilon_i}$.

For such $\epsilon$ and $b_i$, let
$d_i=\delta_i+\epsilon_i$ and choose $\eta_i\in\E_{d_i}$.
The graded tensor product multiplication gives
\[
g(b_1\otimes\eta_1,\cdots,b_n\otimes\eta_n)
=
(-1)^{\sum_{i<j}d_i\epsilon_j}
(\tau_\epsilon h)(b_1,\cdots,b_n)
\otimes\eta_1\cdots\eta_n.
\]
Here the signs follow from
$d_i\epsilon_j+d_j\epsilon_i+d_id_j
=\delta_i\delta_j+\epsilon_i\epsilon_j$ in $\Ztwo$.

Now we show that $g$ is a graded
identity of $B\mathbin{\widehat{\otimes}}\E$ if and only if
$(\tau_\epsilon h)(b_1,\cdots,b_n)=0$ for every
$\epsilon\in(\Ztwo)^n$ and all $b_i\in B_{\epsilon_i}$.
By multilinearity, it suffices to consider
$g(b_1\otimes\eta_1,\cdots,b_n\otimes\eta_n)$,
where $b_i\in B_{\epsilon_i}$ and
$\eta_i\in\E_{\delta_i+\epsilon_i}$ for arbitrary
$\epsilon\in(\Ztwo)^n$.
 Moreover, choosing the odd $\eta_i$
to be distinct exterior generators and the even $\eta_i$
to be $1$ makes $\eta_1\cdots\eta_n$ nonzero.
This proves the lemma.
\end{proof}

The sign transformation also gives the following stability property
of Grassmann envelope identities.

\begin{corollary}\label{cor:envelope-stability}
For every superalgebra $B$ and every $N\geq0$,
\[
\Id\bigl(\E(B\mathbin{\widehat{\otimes}}\E^{(N)})\bigr)=\Id(\E(B)).
\]
\end{corollary}
\begin{proof}
By Theorem~\ref{thm:linearization}, it suffices to compare
multilinear identities. There are graded isomorphisms $
(B\mathbin{\widehat{\otimes}}\E^{(N)})\mathbin{\widehat{\otimes}}\E
\cong B\mathbin{\widehat{\otimes}}(\E^{(N)}\mathbin{\widehat{\otimes}}\E)
\cong B\mathbin{\widehat{\otimes}}\E$.
Applying Lemma~\ref{lem:sign-transfer} with $B$ and
$B\mathbin{\widehat{\otimes}}\E^{(N)}$ therefore gives the same ordinary
multilinear identities for their Grassmann envelopes.
\end{proof}

\subsection{Reduction to graded tensor products
with finite exterior algebras}\label{subsec:finite-exterior-reduction}

We next compare graded identities in a fixed number of variables.
The required number of exterior generators can be bounded using
only the dimension of the semisimple quotient and the nilpotency
index of the radical, independently of the polynomial degree.

\begin{lemma}\label{lem:finite-exterior}
Let $B$ be a finite-dimensional superalgebra over  $\kk$, and let $d\geq1$. Write $s=\dim(B/J(B))$ and choose $\nu\geq1$ such that $J(B)^\nu=0$. If $N\geq ds+\nu-1$, then
$B\mathbin{\widehat{\otimes}}\E^{(N)}$ has the same graded identities as
$B\mathbin{\widehat{\otimes}}\E$ in at most $d$ 
variables.
\end{lemma}
\begin{proof}

Choose a homogeneous basis $b_1,\cdots,b_t$ of $B$ such that
$b_{s+1},\cdots,b_t$ form a basis of $J(B)$.
Let $f=f(x_1,\cdots,x_m)$, with $m\leq d$, be a graded
polynomial.
Write $\delta_i$ for the parity of $x_i$.
For each $i$ and $a$, introduce a formal variable $z_{ia}$
of parity $|z_{ia}|=\delta_i+|b_a|$.
Let $R$ be the tensor product over $\kk$ of the polynomial
algebra in the even variables $z_{ia}$ and the exterior
algebra in the odd variables $z_{ia}$.
Thus 
$z_{ia}z_{jb}=(-1)^{|z_{ia}||z_{jb}|}z_{jb}z_{ia}$,
with $z_{ia}^2=0$ whenever $|z_{ia}|=1$.

For each $i$, define
$X_i=\sum_{a=1}^t b_a\otimes z_{ia}\in B\mathbin{\widehat{\otimes}} R$, which is homogeneous
of parity $\delta_i$.
Given any homogeneous substitution
$x_i\mapsto u_i\in(B\mathbin{\widehat{\otimes}}\E)_{\delta_i}$, write
$u_i=\sum_{a=1}^t b_a\otimes\eta_{ia}$, where
$\eta_{ia}\in\E_{\delta_i+|b_a|}$.
Since $\E$ is supercommutative, the assignments
$z_{ia}\mapsto\eta_{ia}$ define a grading-preserving algebra
homomorphism $\phi:R\to\E$.
The induced homomorphism
$\operatorname{id}_B\otimes\phi:B\mathbin{\widehat{\otimes}} R
\to B\mathbin{\widehat{\otimes}}\E$ sends each $X_i$ to $u_i$.

Suppose $f$ is not a graded identity of $B\mathbin{\widehat{\otimes}}\E$. Then $f(X_1,\cdots,X_m)\ne0$.
 Let $\mathbf y$
denote the even variables $z_{ia}$, and write
\[
f(X_1,\cdots,X_m)
=
\sum_{a,M} b_a\otimes p_{a,M}(\mathbf y)M,
\qquad p_{a,M}(\mathbf y)\in\kk[\mathbf y],
\]
where $M$ ranges over products of distinct odd variables
$z_{ia}$ in lexicographic order of the pairs $(i,a)$, including the empty product $1$.
Since $J(B)^\nu=0$, a product containing at least $\nu$ radical basis elements vanishes.
Thus, for $M=z_{i_1,a_1}\cdots z_{i_q,a_q}$,
$
p_{a,M}\ne0$ implies that $
\#\{j\mid a_j>s\}\leq\nu-1$.
Moreover, the pairs $(i_j,a_j)$ are distinct, and
$1\leq i_j\leq m$, so
$
\#\{j\mid a_j\leq s\}\leq ms\leq ds$.
Hence $q\leq ds+\nu-1$.

Choose $a_0$ and
$M_0=z_{i_1,a_1}\cdots z_{i_q,a_q}$ such that
$p_{a_0,M_0}(\mathbf y)\ne0$.
There are scalars
$\boldsymbol\lambda=(\lambda_{ia})_{|z_{ia}|=0}$ such that
$p_{a_0,M_0}(\boldsymbol\lambda)\ne0$.
Let $\theta_1,\cdots,\theta_N$ be the generators of
$\E^{(N)}$.
Since $q\leq ds+\nu-1\leq N$, the assignments
\[
\psi(z_{ia})=
\begin{cases}
\lambda_{ia},
  & |z_{ia}|=0,\\
\theta_j,
  & (i,a)=(i_j,a_j),\quad 1\leq j\leq q,\\
0,
  & |z_{ia}|=1,\quad
    (i,a)\notin\{(i_1,a_1),\cdots,(i_q,a_q)\}
\end{cases}
\]
define a grading-preserving algebra homomorphism
$\psi:R\to\E^{(N)}$.
Let $v_i=(\operatorname{id}_B\otimes\psi)(X_i)$.
Then
$
f(v_1,\cdots,v_m)
=
(\operatorname{id}_B\otimes\psi)
(f(X_1,\cdots,X_m))$.
In the tensor product basis of $B\mathbin{\widehat{\otimes}}\E^{(N)}$,
the coefficient of
$b_{a_0}\otimes \prod_{j=1}^q\theta_j$
in this value is
$p_{a_0,M_0}(\boldsymbol\lambda)\ne0$.
Thus $f(v_1,\cdots,v_m)\ne0$, so $f$ is not a graded
identity of $B\mathbin{\widehat{\otimes}}\E^{(N)}$.
The converse follows from $\E^{(N)}\subseteq\E$.

\end{proof}

We also need a bound on the nilpotency index of
$J(B\mathbin{\widehat{\otimes}}\E^{(N)})$.

\begin{lemma}\label{lem:tensor-radical}
Let $B$ be a finite-dimensional superalgebra, and let $N\geq0$ and $\nu\geq1$ be integers. If $J(B)^\nu=0$, then
$
J(B\mathbin{\widehat{\otimes}}\E^{(N)})^{\nu+N}=0$.
More precisely, writing $\E^{(N)}_+$ for the ideal spanned by
exterior products of positive degree, this radical is
$
J(B)\mathbin{\widehat{\otimes}}\E^{(N)}
  +B\mathbin{\widehat{\otimes}}\E^{(N)}_+$.
\end{lemma}

\begin{proof}
Let $C=B\mathbin{\widehat{\otimes}}\E^{(N)}$ and $
R=
J(B)\mathbin{\widehat{\otimes}}\E^{(N)}
+B\mathbin{\widehat{\otimes}}\E^{(N)}_+$.
Let $\varepsilon:\E^{(N)}\to\kk$ send each element to its
constant term. The algebra homomorphism $
\pi:C\longrightarrow B/J(B)$ defined by 
$
\pi(b\otimes\eta)=\varepsilon(\eta)(b+J(B))$ 
is surjective, and 
its kernel is
$J(B)\mathbin{\widehat{\otimes}}\kk1
+B\mathbin{\widehat{\otimes}}\E^{(N)}_+=R$.
Thus $R$ is a two-sided ideal and $C/R\cong B/J(B)$. Since $J(B)^\nu=0$, it is easy to check that  $R^{\nu+N}=0$.

Since $R$ is a nilpotent ideal, $R\subseteq J(C)$.
On the other hand, $C/R\cong B/J(B)$ is semisimple,
so $J(C)\subseteq R$.
Hence $J(C)=R$, and consequently $J(C)^{\nu+N}=0$.
\end{proof}

% ===== parts/separation =====
\section{Separating polynomials}\label{sec:separation}

Throughout this section, $\kk$ is algebraically closed. Let
$\mathcal S=(S_1,\cdots,S_r)$ be a nonempty ordered family of
finite-dimensional graded-simple superalgebras, and let $I(\mathcal S)$
be as in Definition~\ref{def:ordered-family}.
We construct graded polynomials $F_L$ in a number of variables
independent of $L$ and give nonzero graded evaluations on a minimal
superalgebra associated with $\mathcal S$.
Conversely, for sufficiently large $L$, nonzero graded evaluations
of $F_L$ on finite-dimensional superalgebras yield minimal
superalgebras associated with $\mathcal S$.
These constructions will be used to characterize containment
in $I(\mathcal S)$.

\subsection{The construction of separating polynomials}\label{subsec:separating-construction}

We first define $F_L$ and the associated ordinary multilinear
polynomial $f_{\mathcal S,L}$. We then construct a minimal superalgebra
on which every $F_L$ has a nonzero graded evaluation and prove that
$f_{\mathcal S,L}\notin I(\mathcal S)$.

\begin{definition}\label{def:separator-data}
For $S_i=\Mat_{p_i|q_i}(\kk)$, let $a_i=p_i+q_i$.
Assign parities to the matrix indices, with $p_i$ even indices
and $q_i$ odd indices. Let $c_{ij}\in\Ztwo$ denote the parity
of index $j$.
The matrix unit $E_{jk}$ has parity $c_{ij}+c_{ik}$.
Replacing every $c_{ij}$ by $c_{ij}+1$ leaves these
sums unchanged, so
$\Mat_{p_i|q_i}(\kk)\cong\Mat_{q_i|p_i}(\kk)$
as graded algebras.
After interchanging $p_i$ and $q_i$ if necessary,
we may order the indices so that $c_{i1}=0$.
Set $c_{i,a_i+1}=c_{i1}$ and introduce variables
$x_{i1},\cdots,x_{i,a_i}$ with
$
|x_{ij}|=c_{ij}+c_{i,j+1}$
for $1\leq j\leq a_i$.

For $S_i=\Qtype_{n_i}$, introduce variables
$x^+_{ij},x^-_{ij}$ $(1\leq j\leq n_i)$ and
$\xi_i,\eta_i$, with
$
|x^+_{ij}|=|x^-_{ij}|=0$ and $
|\xi_i|=|\eta_i|=1$.

For $S_i=\Mat_{p_i|q_i}(\kk)$, let
$\mathcal R_i=\{(i,j)\mid 1\leq j\leq a_i\}$.
For $S_i=\Qtype_{n_i}$, let
\[
\mathcal R_i^+=\{(i,+,j)\mid 1\leq j\leq n_i\},
\qquad
\mathcal R_i^-=\{(i,-,j)\mid 1\leq j\leq n_i\},
\qquad
\mathcal R_i=\mathcal R_i^+\cup\mathcal R_i^-.
\]
Let $\mathcal R=\bigcup_{i=1}^r\mathcal R_i$.
The labels $(i,j)$, $(i,+,j)$, and $(i,-,j)$ denote
the positions immediately before $x_{ij}$, $x^+_{ij}$,
and $x^-_{ij}$, respectively.
Order $\mathcal R$ first by increasing $i$.
Explicitly, the order within $\mathcal R_i$ is $
(i,1)<\cdots<(i,a_i)$
for $S_i=\Mat_{p_i|q_i}(\kk)$, and
$
(i,+,1)<\cdots<(i,+,n_i)
<(i,-,1)<\cdots<(i,-,n_i)
$
for $S_i=\Qtype_{n_i}$.

Define
\[
\mathcal P
=
\{(a,b)\in\mathcal R\times\mathcal R\mid a<b\}
\setminus
\left(
\bigcup_{S_i=\Qtype_{n_i}}
(\mathcal R_i^+\times\mathcal R_i^-)
\right).
\]
Introduce commuting indeterminates $Z_a$ for $a\in\mathcal R$,
and define
$
H(Z)=\prod_{(a,b)\in\mathcal P}(Z_b-Z_a)$.
Write
$
H(Z)=\sum_{\mathbf k}
c_{\mathbf k}\prod_{a\in\mathcal R}Z_a^{k_a}$,
where $\mathbf k=(k_a)_{a\in\mathcal R}$ ranges over
$\mathbb Z_{\geq0}^{\mathcal R}$ and $c_{\mathbf k}\in\mathbb Z$.
Only finitely many $c_{\mathbf k}$ are nonzero.
Since $H$ is homogeneous, $c_{\mathbf k}\ne0$ implies
$\sum_{a\in\mathcal R}k_a=\deg H$.

\end{definition}

\begin{definition}\label{def:separating-polynomials}
Using the data of Definition~\ref{def:separator-data}, introduce
an even variable $t$.
For each $\mathbf k\in\mathbb Z_{\geq0}^{\mathcal R}$, define
\[
W_i(\mathbf k)
=
\prod_{j=1}^{a_i}
\bigl(t^{k_{(i,j)}}x_{ij}\bigr)
\qquad
\text{if }S_i=\Mat_{p_i|q_i}(\kk),
\]
and
\[
W_i(\mathbf k)
=
\left(
\prod_{j=1}^{n_i}
\bigl(t^{k_{(i,+,j)}}x^+_{ij}\bigr)
\right)\xi_i
\left(
\prod_{j=1}^{n_i}
\bigl(t^{k_{(i,-,j)}}x^-_{ij}\bigr)
\right)\eta_i
\qquad
\text{if }S_i=\Qtype_{n_i}.
\]
The products are taken in increasing order of $j$.
Whenever $k_a=0$, the corresponding power of $t$ is omitted.

Introduce even variables $y_1,\cdots,y_{r-1}$.
For $L\geq1$, define
\begin{equation}\label{eq:FL}
F_L=
\sum_{\mathbf k^{(1)},\cdots,\mathbf k^{(L)}}
\left(\prod_{\ell=1}^L c_{\mathbf k^{(\ell)}}\right)
\left[
\prod_{i=1}^{r-1}
\left(
\left(\prod_{\ell=1}^L W_i(\mathbf k^{(\ell)})\right)y_i
\right)
\right]
\left(\prod_{\ell=1}^L W_r(\mathbf k^{(\ell)})\right).
\end{equation}
Each $\mathbf k^{(\ell)}$ ranges over
$\mathbb Z_{\geq0}^{\mathcal R}$.
Only finitely many summands are nonzero.
The products indexed by $i$ and $\ell$ are taken in
increasing order. For $r=1$, the product indexed by $i$
is omitted.
For each $\ell$, the factors
$W_1(\mathbf k^{(\ell)}),\cdots,W_r(\mathbf k^{(\ell)})$
form layer $\ell$.
All factors use the variables from Definition~\ref{def:separator-data}
and the same variable $t$.

The polynomial $F_L$ is multihomogeneous.
It has degree $L$ in each variable $x_{ij}$, $x^+_{ij}$,
$x^-_{ij}$, $\xi_i$, or $\eta_i$ that occurs, degree one
in each $y_i$, and degree $L\deg H$ in $t$.
The number of variables is bounded independently of $L$ by
\begin{equation}\label{eq:variable-bound}
d=r+\sum_{S_i=\Mat_{p_i|q_i}}(p_i+q_i)
  +\sum_{S_i=\Qtype_{n_i}}(2n_i+2).
\end{equation}
The term $r$ counts $t$ and the $r-1$ variables $y_i$.
If $H=1$, then $t$ is absent, and $d$ remains an upper bound.

Let $g_L$ be the complete linearization of $F_L$,
with each new variable having the same parity as
the variable it replaces.
Let $\delta_L$ be the parity vector of the variables in $g_L$.
Define the ordinary multilinear polynomial
\[
f_{\mathcal S,L}=\tau_{\delta_L}g_L.
\]
The coefficients of $F_L$, $g_L$, and $f_{\mathcal S,L}$
are integers.
\end{definition}

%Different from the strongly full polynomials used by
%Aljadeff and Karasik \cite{AljadeffKarasik2024},
%which are constructed using alternation and
%symmetrization, our polynomials $F_L$ are based
%on a Vandermonde-type construction and involve
%a number of variables independent of $L$.

We realize $\mathcal S$ as the diagonal blocks of a superalgebra
on which each $F_L$ has a nonzero graded evaluation.

\begin{lemma}\label{lem:target-evaluation}
Let
\[
m_i=
\begin{cases}
p_i+q_i, & S_i=\Mat_{p_i|q_i}(\kk),\\
2n_i,   & S_i=\Qtype_{n_i},
\end{cases}
\]
and let
\[
A(\mathcal S)=
\begin{pmatrix}
S_1 & \Mat_{m_1\times m_2}(\kk) & \cdots
    & \Mat_{m_1\times m_r}(\kk)\\
0   & S_2 & \cdots & \Mat_{m_2\times m_r}(\kk)\\
\vdots & \ddots & \ddots & \vdots\\
0 & \cdots & 0 & S_r
\end{pmatrix},
\]
where each $S_i$ is regarded as a subalgebra of
$\Mat_{m_i}(\kk)$ via the embeddings described earlier.
Then $A(\mathcal S)$ is minimal with respect to
$(S_1,\cdots,S_r)$, and every $F_L$ has a nonzero graded
evaluation on $A(\mathcal S)$.
 Moreover, $
\Id(\E(A(\mathcal S)))=I(\mathcal S)$, and $ f_{\mathcal S,L}\notin I(\mathcal S)$ for $L\geq1$.

\end{lemma}

\begin{proof}

Let $A=A(\mathcal S)$, and let $E_{ab}$ denote the standard
matrix units of $\Mat_{m_1+\cdots+m_r}(\kk)$. Write $c_a$ for the
parity of index $a$. Then $A=A_0\oplus A_1$, where
$
A_0=A\cap\operatorname{span}_\kk\{E_{ab}\mid c_a=c_b\}$ and $
A_1=A\cap\operatorname{span}_\kk\{E_{ab}\mid c_a\ne c_b\}$.
Let $R$ be the subspace of strictly upper block matrices.
Then $R$ is a graded ideal, $R^r=0$, and $
A/R\cong S_1\oplus\cdots\oplus S_r$.
Since each $S_i$ is semisimple as an ungraded algebra,
we have $R=J(A)$.

For a matrix-type component $S_i$, let $E_{ab}^{(i)}$
denote its standard matrix units.
For $S_i=\Qtype_{n_i}$, let $E_{ab}^{(i)}$ denote the
standard matrix units of its even subalgebra
$\Mat_{n_i}(\kk)$, and write $u_i$ for the odd element
satisfying $u_i^2=1$ and commuting with this subalgebra.
In both cases, let $e_i=E_{11}^{(i)}$.
Thus $e_i$ is a minimal nonzero homogeneous idempotent
of $S_i$.
In the matrix embedding of $\Qtype_{n_i}$, it has
diagonal entries equal to one in positions $1$ and
$n_i+1$, and zero elsewhere.

Suppose first that $r\geq2$.
For $1\leq i<r$, let $w_i$ be the matrix whose only
nonzero entry is one in the first row and first column
of the $(i,i+1)$ block.
The first matrix index of every diagonal block is even,
so $w_i$ is even.
Moreover, $
w_i=e_iw_ie_{i+1}\in e_iRe_{i+1}$, and $
w_1\cdots w_{r-1}\ne0$.
It is easy to check that $R$ is generated by
$w_1,\cdots,w_{r-1}$ as an ideal of $A$.
Together with $A=(S_1\oplus\cdots\oplus S_r)\oplus R$,
this proves that $A$ is minimal with respect to $\mathcal S$.
For $r=1$, minimality follows from $A=S_1$.

We next construct a graded evaluation that is nonzero
for every $F_L$.
For $S_i=\Mat_{p_i|q_i}(\kk)$, substitute
$x_{ij}=E_{j,j+1}^{(i)}$ for $1\leq j<a_i$
and $x_{i,a_i}=E_{a_i,1}^{(i)}$.
When $a_i=1$, the substitution is simply
$x_{i1}=E_{11}^{(i)}$.
For $S_i=\Qtype_{n_i}$, substitute
$x^+_{ij}=x^-_{ij}=E_{j,j+1}^{(i)}$ for $1\leq j<n_i$
and $x^+_{i,n_i}=x^-_{i,n_i}=E_{n_i,1}^{(i)}$.
Also substitute $\xi_i=\eta_i=u_ie_i$. Each variable is assigned a homogeneous element
of the same parity.
Under these substitutions, $W_i(\mathbf0)$ evaluates
to $e_i$ in both cases.

Choose scalars $\lambda_{ij}$, with
$1\leq j\leq a_i$ for a matrix-type component and
$1\leq j\leq n_i$ for a $\Qtype$-component,
such that all the chosen scalars are pairwise distinct.
Let
\[
D_i=
\begin{cases}
\displaystyle\sum_{j=1}^{a_i}\lambda_{ij}E_{jj}^{(i)},
& S_i=\Mat_{p_i|q_i}(\kk),\\[6pt]
\displaystyle\sum_{j=1}^{n_i}\lambda_{ij}E_{jj}^{(i)},
& S_i=\Qtype_{n_i},
\end{cases}
\]
and substitute
$t=\operatorname{diag}(D_1,\cdots,D_r)$.
This is an even element of $A$.
For the labels in $\mathcal R$, set $
\lambda_{(i,j)}=\lambda_{ij} $ if $S_i=\Mat_{p_i|q_i}(\kk)$,
and $
\lambda_{(i,+,j)}=\lambda_{(i,-,j)}=\lambda_{ij}$ if $S_i=\Qtype_{n_i}$.
Under these substitutions, each inserted power $t^{k_a}$
multiplies the following matrix unit by $\lambda_a^{k_a}$.
Since $W_i(\mathbf0)$ evaluates to $e_i$, the value of
$W_i(\mathbf k)$ is $(\prod_{a\in\mathcal R_i}\lambda_a^{k_a}) e_i$.

The only repeated entries in
$\lambda=(\lambda_a)_{a\in\mathcal R}$ occur between
the two segments of the same $\Qtype$-component.
The corresponding pairs are excluded from $\mathcal P$.
Consequently,
$
H(\lambda)=\prod_{(a,b)\in\mathcal P}
(\lambda_b-\lambda_a)\ne0$.
Finally, substitute $y_i=w_i$.
Substituting these values of $W_i(\mathbf k)$ in
\eqref{eq:FL}, 
we have that the value of $F_L$ is $H(\lambda)^L e_1w_1e_2\cdots w_{r-1}e_r\ne0$.
For $r=1$, the value is $H(\lambda)^L e_1$.
Thus every $F_L$ has a nonzero graded evaluation on $A$.

By minimality and Theorem~\ref{thm:minimal-product},
$
\Id(\E(A))
=\Id(\E(S_1))\cdots\Id(\E(S_r))
=I(\mathcal S)$.
By Theorem~\ref{thm:linearization}, $g_L$ is not
a graded identity of $A$.
The graded embedding $
A\hookrightarrow A\mathbin{\widehat{\otimes}}\E$ defined by $ a\longmapsto a\otimes1$
shows that $g_L$ is not a graded identity of
$A\mathbin{\widehat{\otimes}}\E$.
Lemma~\ref{lem:sign-transfer} therefore gives $
f_{\mathcal S,L}=\tau_{\delta_L}g_L
\notin\Id(\E(A))=I(\mathcal S)$.

\end{proof}

\begin{remark}\label{rem:target-evaluation-field}
Let $\kk$ be any field of characteristic zero, and let
$\mathcal S$ be a nonempty ordered family of split graded-simple
superalgebras over $\kk$.
The constructions of $F_L$ and $A(\mathcal S)$ are defined
over $\mathbb Q$, and the distinct scalars $\lambda_{ij}$
can be chosen in $\mathbb Q$.
Thus every $F_L$ has a nonzero graded evaluation on $A(\mathcal S)$.
Over an algebraic closure $\overline\kk$, the superalgebra
$A(\mathcal S)\otimes_\kk\overline\kk$ is minimal with respect to
$(S_1\otimes_\kk\overline\kk,\cdots,S_r\otimes_\kk\overline\kk)$
by Lemma~\ref{lem:target-evaluation}.
Theorem~\ref{thm:minimal-product} and Lemma~\ref{lem:scalar-extension}
therefore give $\Id(\E(A(\mathcal S)))=I(\mathcal S)$.
Theorem~\ref{thm:linearization} and Lemma~\ref{lem:sign-transfer}
then give $f_{\mathcal S,L}\notin I(\mathcal S)$ for every $L\geq1$.
\end{remark}

\subsection{Minimal superalgebras arising from nonzero evaluations}\label{subsec:nonzero-evaluations}

We now recover a minimal superalgebra from a nonzero graded evaluation
of $F_L$ for sufficiently large $L$. 
The first step is to put the value
substituted for $t$ in a normal form relative to a graded
Wedderburn--Malcev decomposition.

\begin{lemma}\label{lem:diagonal-radical}
Let $B=S\oplus J$ be a graded Wedderburn--Malcev decomposition
of a finite-dimensional superalgebra, where $J=J(B)$ and
$J^\nu=0$.
Let $f_a$ denote the diagonal matrix units of $S_0$.
Let $U\subseteq B_0$ consist of the elements whose projections
to $S_0$ have nonzero eigenvalues that are pairwise distinct
among all matrix summands.
Then $U$ is Zariski dense and open.
For every $T\in U$, there is an even invertible element $g$
in the unitization of $B$ such that
\[
gTg^{-1}=D+\zeta,\qquad
D=\sum_a\lambda_a f_a,\qquad
\zeta\in J_0,\qquad [D,\zeta]=0,
\]
where the $\lambda_a$ are pairwise distinct and nonzero.
If $B$ is nonunital, set $f_0=1-1_S$ and $\lambda_0=0$
in the unitization.
Then $\zeta$ commutes with every $f_a$, including $f_0$
when present.
\end{lemma}

\begin{proof}
The eigenvalue conditions define a nonempty Zariski open
subset of $S_0$. Its inverse image $U$ is therefore open
and dense in $B_0$.

We work in the unitization of $B$, understood to be $B$
itself when $B$ is already unital, and denote its identity
by $1$.
If $B$ is unital, then $1_S=1$. Indeed, the image of $1_S$
in $B/J$ is its identity, so $1-1_S$ is an idempotent
in the nilpotent ideal $J$ and must be zero.
If $B$ is nonunital, we extend the family $(f_a)$ by
the idempotent $f_0=1-1_S$ and set $\lambda_0=0$,
as in the statement. In this case, all sums below
are taken over the extended index set.
Thus, in either case, $\sum_a f_a=1$.

For $T\in U$, diagonalizing its projection to $S_0$
gives an even invertible element $h$ such that
$hTh^{-1}=\sum_a\lambda_a f_a+z=D+z$, with $z\in J_0$.
In the nonunital case, extend the conjugating element
from $S_0$ by adding $f_0$.

We prove by induction on $m$ that $D+z$ is conjugate,
by an even invertible element, to $D+z_m$, where
$z_m\in J_0$ and $f_a z_m f_b\in J^m$ for all $a\ne b$.
For $m=1$, take $z_1=z$.
Suppose the assertion holds for $m<\nu$, and define
\[
v_m=\sum_{a\ne b}
\frac{f_a z_m f_b}{\lambda_a-\lambda_b}
\in(J^m)_0.
\]
Then $1+v_m$ is even and invertible. Let
$z_{m+1}=(1+v_m)(D+z_m)(1+v_m)^{-1}-D$. Since $1+v_m\equiv1\pmod J$ and all these elements are even,
we have $z_{m+1}\in J_0$.
 Since $v_m\in J^m$ and $z_m\in J$,
$
z_{m+1}\equiv z_m+[v_m,D]
=\sum_a f_a z_m f_a
\pmod{J^{m+1}}$.
Thus $f_a z_{m+1}f_b\in J^{m+1}$ for $a\ne b$,
completing the induction.

Taking $\zeta=z_\nu$ and using $J^\nu=0$, we obtain
$\zeta=\sum_a f_a\zeta f_a$.
Hence $\zeta$ commutes with every $f_a$ and with $D$.
Composing these conjugations with the one by $h$
gives the required element $g$.

\end{proof}

We apply this normal form to a nonzero graded evaluation of $F_L$
and use its expansion to construct the required minimal superalgebra.

\begin{lemma}\label{lem:extraction}
Let
$B=(R_1\oplus\cdots\oplus R_q)\oplus J$
be a graded Wedderburn--Malcev decomposition of a finite-dimensional
superalgebra over an algebraically closed field,
where the $R_j$ are graded-simple, $J=J(B)$, and $J^\nu=0$.
Suppose that $L\geq\nu$ and $F_L$ is not a graded identity
of $B$.

Let $\mathcal S=(S_1,\cdots,S_r)$ be the ordered family
of graded-simple superalgebras in
Definition~\ref{def:separating-polynomials}. Then there are indices $\alpha_i\in\{1,\cdots,q\}$
and graded subalgebras $T_i\subseteq R_{\alpha_i}$
such that $T_i\cong S_i$ and the identity elements
of the $T_i$ are pairwise orthogonal.
The indices $\alpha_i$ need not be distinct.

There are minimal homogeneous idempotents $e_i\in T_i$
and homogeneous elements $v_i\in e_iBe_{i+1}$
such that $v_1\cdots v_{r-1}\ne0$.
The subalgebras $T_i$ and the elements $v_i$ generate
a minimal superalgebra $C\subseteq B$. For $r=1$, take $C=T_1$ and omit the elements $v_i$. In particular, $
\Id(\E(B))\subseteq I(\mathcal S)$.

\end{lemma}

\begin{proof}
Let $d_0$ and $d_1$ be the numbers of distinct even and odd
variables in $F_L$, respectively.
Since $F_L$ is not a graded identity of $B$, the set
$\{\mathbf b\in B_0^{d_0}\times B_1^{d_1}\mid
F_L(\mathbf b)\ne0\}$
is a nonempty Zariski open subset of
$B_0^{d_0}\times B_1^{d_1}$.

Let $U\subseteq B_0$ be the dense open subset in
Lemma~\ref{lem:diagonal-radical}.
If $t$ occurs in $F_L$, the tuples for which the value
assigned to $t$ lies in $U$ form a dense open subset of
$B_0^{d_0}\times B_1^{d_1}$.
Thus we may choose such a tuple $\mathbf b$ with
$F_L(\mathbf b)\ne0$.
If $t$ does not occur in $F_L$, choose $\mathbf b$ with
$F_L(\mathbf b)\ne0$ and choose any $T\in U$
as the value of $t$.

Let $f_b$ be the diagonal matrix units of
$(R_1\oplus\cdots\oplus R_q)_0$.
By Lemma~\ref{lem:diagonal-radical}, after conjugating
all chosen elements by the same even invertible element $g$,
we may assume that the value $T$ assigned to $t$
has the form
\[
T=D+\zeta,\qquad
D=\sum_b\lambda_b f_b,\qquad
\zeta\in J_0,\qquad f_b\zeta=\zeta f_b,
\]
where the $\lambda_b$ are pairwise distinct and nonzero.
The value of $F_L$ is conjugated by $g$
and therefore remains nonzero.
If $B$ is nonunital, include $f_0$ and $\lambda_0=0$
as in Lemma~\ref{lem:diagonal-radical}.

Choose homogeneous bases of the spaces $f_bJf_c$.
In each $R_j$ of matrix type, use the homogeneous
matrix units.
If $R_j=\Qtype_m$, use the elements $E_{ab}$ and $uE_{ab}$.
Together these form a homogeneous basis of $B$.
Express the value assigned to each variable other than $t$
as a linear combination of these basis elements,
and expand each product by distributivity.
At different occurrences of the same variable, basis
elements are chosen independently from that
linear combination, and their coefficients are multiplied.

For a chosen basis element $x$ with $f_bx=x$, we have
$Dx=\lambda_bx$.
Since $D\zeta=\zeta D$ and $J^\nu=0$, the binomial formula gives
\begin{equation}\label{eq:binomial-radical}
T^kx
=\sum_{h=0}^{\min(k,\nu-1)}
\binom{k}{h}\lambda_b^{k-h}\zeta^h x.
\end{equation}
Here $\lambda_b^0=1$, also when $\lambda_b=0$.

Choose one basis element from the linear combination at each
occurrence of a variable other than $t$, and one exponent $h$
in each sum in \eqref{eq:binomial-radical}.
We now compute the coefficient corresponding to these choices.
For each layer $\ell$ and recorded position $a$,
let $h_a^{(\ell)}$ be the chosen exponent, and let
$\lambda_a^{(\ell)}=\lambda_b$ if the chosen basis element
$x$ at that position satisfies $f_bx=x$.
Write
$\mathbf h^{(\ell)}=(h_a^{(\ell)})_{a\in\mathcal R}$
and
$\lambda^{(\ell)}=(\lambda_a^{(\ell)})_{a\in\mathcal R}$.
For these choices, let $W$ be the product of
the selected basis elements and powers of $\zeta$,
in their original order, and let $\gamma$ be the
product of the coefficients of the selected basis
elements, one for each occurrence of a variable
other than $t$.
Neither $W$ nor $\gamma$ depends on the exponent
vectors $\mathbf k^{(\ell)}$.
Summing over the exponent vector $\mathbf k^{(\ell)}$
in \eqref{eq:FL} gives the scalar
$H_{\mathbf h^{(\ell)}}(\lambda^{(\ell)})$, where
\[
H_{\mathbf h}(\lambda)
=
\left(\prod_{a\in\mathcal R}\frac{1}{h_a!}\right)
\left(\prod_{a\in\mathcal R}\partial_{Z_a}^{h_a}\right)
H(Z)\bigg|_{Z_a=\lambda_a}.
\]
Indeed, for $k_a\geq h_a$, differentiating $Z_a^{k_a}$
and dividing by $h_a!$ gives
$\binom{k_a}{h_a}Z_a^{k_a-h_a}$.
After substituting $Z_a=\lambda_a$, this is precisely
the coefficient in \eqref{eq:binomial-radical}.
For $k_a<h_a$, the derivative is zero.

Thus the value of $F_L$ is a finite sum of terms
\[
\gamma
\left(
\prod_{\ell=1}^L
H_{\mathbf h^{(\ell)}}(\lambda^{(\ell)})
\right)W.
\]
Choose a nonzero term in this sum.
Then $W\ne0$ and
$H_{\mathbf h^{(\ell)}}(\lambda^{(\ell)})\ne0$
for every $\ell$.

Let $c$ be the number of occurrences at which a basis
element of $J$ was chosen, including those among the
variables $y_i$.
Since $W\ne0$, $\zeta\in J$, and $J^\nu=0$, we have $
c+\sum_{\ell=1}^L\sum_{a\in\mathcal R}
h_a^{(\ell)}<\nu$.
Each layer containing a chosen basis element of $J$
or a positive exponent $h_a^{(\ell)}$ contributes
at least one to this count.
The variables $y_i$ belong to no layer.
Hence there are at most $\nu-1$ such layers.
Since $L\geq\nu$, there is a layer $\ell_0$ in which
all chosen basis elements belong to
$R_1\oplus\cdots\oplus R_q$
and all $h_a^{(\ell_0)}$ are zero.
For this layer, $H(\lambda^{(\ell_0)})\ne0$.

For each $i$, let $w_i$ be the part of $W$ corresponding to
$W_i(\mathbf k^{(\ell_0)})$.
Since all $h_a^{(\ell_0)}$ are zero, $w_i$ is a product
of the chosen basis elements in
$R_1\oplus\cdots\oplus R_q$.
Moreover, $w_i\ne0$, since it is a consecutive part
of the nonzero product $W$.
As elements of distinct components $R_j$ multiply to zero,
all basis elements in $w_i$ belong to one component
$R_{\alpha_i}$.
We now construct a graded subalgebra
$T_i\subseteq R_{\alpha_i}$ with $T_i\cong S_i$.

Suppose first that $S_i=\Mat_{p_i|q_i}$, and let
$a_i=p_i+q_i$.
Let $\rho_1,\cdots,\rho_{a_i}$ be the row indices
of the basis elements chosen for
$x_{i1},\cdots,x_{i,a_i}$ in $w_i$.
These indices are pairwise distinct because
$H(\lambda^{(\ell_0)})\ne0$.

If $R_{\alpha_i}$ is of matrix type, the nonvanishing
of $w_i$ implies that the column index of each chosen
matrix unit agrees with the row index of the next.
The given parities of the variables $x_{ij}$
therefore imply that the parities of
$\rho_1,\cdots,\rho_{a_i}$
are $c_{i1},\cdots,c_{i,a_i}$, up to reversing all parities.
Hence the span of the matrix units $E_{\rho_a,\rho_b}$
is a graded subalgebra $T_i\cong\Mat_{p_i|q_i}$.

If $R_{\alpha_i}=\Qtype_m$, define a map by
$
E_{ab}\longmapsto
u^{c_{ia}+c_{ib}}E_{\rho_a,\rho_b}$.
This map preserves degrees and is injective because
the $\rho_a$ are distinct.
It preserves multiplication because $u$ commutes with
the matrix units and $u^2=1$.
Its image is therefore a graded subalgebra
$T_i\cong\Mat_{p_i|q_i}$.

Suppose now that $S_i=\Qtype_{n_i}$.
Let $\rho_1,\cdots,\rho_{n_i}$ and
$\sigma_1,\cdots,\sigma_{n_i}$ be the row indices
of the basis elements chosen for the variables
$x^+_{ij}$ and $x^-_{ij}$, respectively, in $w_i$. Since 
$H(\lambda^{(\ell_0)})\ne0$, the indices $\rho_1,\cdots,\rho_{n_i}$ are pairwise distinct,
as are $\sigma_1,\cdots,\sigma_{n_i}$.

If $R_{\alpha_i}=\Qtype_m$, the span of
$E_{\rho_a,\rho_b}$ and $uE_{\rho_a,\rho_b}$
is a graded subalgebra $T_i\cong\Qtype_{n_i}$.

If $R_{\alpha_i}$ is of matrix type, the indices
$\rho_1,\cdots,\rho_{n_i}$ all have the same parity,
since the variables $x^+_{ij}$ are even and $w_i\ne0$.
The same holds for $\sigma_1,\cdots,\sigma_{n_i}$.
The odd variable $\xi_i$ between the two families
of variables forces these two parities to be opposite.
Thus the sets $\{\rho_1,\cdots,\rho_{n_i}\}$
and $\{\sigma_1,\cdots,\sigma_{n_i}\}$ are disjoint.
Define $
A_{ab}=E_{\rho_a,\rho_b}+E_{\sigma_a,\sigma_b}$ and $
B_{ab}=E_{\rho_a,\sigma_b}+E_{\sigma_a,\rho_b}$.
The elements $A_{ab}$ are even, the elements $B_{ab}$
are odd, and $
A_{ab}A_{cd}=\delta_{bc}A_{ad}$, $
A_{ab}B_{cd}=\delta_{bc}B_{ad}$,  $ 
B_{ab}A_{cd}=\delta_{bc}B_{ad}$,
$B_{ab}B_{cd}=\delta_{bc}A_{ad}$.
Hence their span is a graded subalgebra
$T_i\cong\Qtype_{n_i}$.
These constructions use the row indices at the recorded
positions; the last column index of $w_i$ may differ
from its first row index.

Every pair of recorded positions belonging to different
$W_i$ belongs to $\mathcal P$. Since $H(\lambda^{(\ell_0)})\ne0$, the row indices
at such positions cannot be equal within the same
component $R_j$.
Thus, within each $R_j$, the sets of matrix indices
used for different $T_i$ are disjoint.
It follows that their identity elements
$\epsilon_i=1_{T_i}$ are pairwise orthogonal.

Let $f_{b_i}$ be the diagonal matrix unit corresponding
to the row index of the first basis element of $w_i$.
Define $e_i=f_{b_i}$, except when $T_i$ is an algebra
of type $\Qtype$ constructed in a component of matrix type,
in which case define $e_i=A_{11}$.
Then $e_i$ is a minimal nonzero homogeneous idempotent
of $T_i$, and
$
e_if_{b_i}=f_{b_i}=f_{b_i}e_i$.

For $r\geq2$ and $1\leq i<r$, let $v_i$ be the product
in $W$ starting with the first basis element of $w_i$
and ending immediately before the first basis element
of $w_{i+1}$.
Since $\zeta$ commutes with every $f_b$ and $W\ne0$,
any two consecutive chosen basis elements
$x\in f_bBf_c$ and $y\in f_dBf_e$ must have $c=d$;
otherwise $x\zeta^h y=0$ for every $h\geq0$.
By the definition of $v_i$, this gives
$
f_{b_i}v_i=v_i=v_if_{b_{i+1}}$.
Since $e_if_{b_i}=f_{b_i}$ and
$f_{b_{i+1}}e_{i+1}=f_{b_{i+1}}$, we obtain
$
v_i\in f_{b_i}Bf_{b_{i+1}}
\subseteq e_iBe_{i+1}$.
Every $v_i$ is homogeneous, since the chosen basis elements
are homogeneous and $\zeta$ is even.
Moreover, $v_1\cdots v_{r-1}$ is a consecutive part
of the nonzero product $W$, so
$v_1\cdots v_{r-1}\ne0$.

Let $C$ be the graded subalgebra generated by the $T_i$
and the $v_i$, and let $N$ be the ideal of $C$
generated by the $v_i$. 
Since
$T_i\subseteq\epsilon_iB\epsilon_i$
and
$v_i\in\epsilon_iB\epsilon_{i+1}$,
orthogonality of the $\epsilon_i$ gives $
N\subseteq\bigoplus_{i<j}\epsilon_iB\epsilon_j$,
$N^r=0$,
and $
C=(T_1\oplus\cdots\oplus T_r)\oplus N$.
Thus $N$ is nilpotent and $C/N$ is semisimple,
so $N=J(C)$.
It follows that $C$ is minimal with respect to
$(T_1,\cdots,T_r)$.
The elements $v_i$ belong to $J(C)$, but need not
belong to $J$.
For $r=1$, take $C=T_1$ directly.

Since $C\subseteq B$, Theorem~\ref{thm:minimal-product}
gives
$
\Id(\E(B))
\subseteq\Id(\E(C))
=\prod_{i=1}^r\Id(\E(T_i))
=\prod_{i=1}^r\Id(\E(S_i))
=I(\mathcal S)$.
\end{proof}

\begin{remark}
The proof of Lemma~\ref{lem:extraction} constructs the $T_i$
from pairwise disjoint sets of matrix indices.
It selects positions in the nonzero product $W$, in the order
$1,\cdots,r$, whose row indices belong to the respective
index sets of $T_1,\cdots,T_r$.

\end{remark}

% ===== parts/intersections =====
\section{Intersection theorems for \texorpdfstring{$T$}{T}-ideals}\label{sec:finite-intersections}

Throughout this section, $\kk$ is an arbitrary field of
characteristic zero. We combine Lemma~\ref{lem:finite-exterior}
and Lemma~\ref{lem:extraction} to determine whether the identity ideal
of a Grassmann envelope is contained in $I(\mathcal S)$.
The polynomial $f_{\mathcal S,L}$ gives this criterion with a
choice of $L$ depending only on bounds for the dimension of the
semisimple quotient and the nilpotency index of the radical.
This uniform choice then gives an intersection theorem for
families of identity ideals. 
For finite families, it proves Berele's Conjecture~3,
from which his Conjecture~1 follows \cite{Berele2012}.

\begin{theorem}\label{thm:uniform-separation}
Let $\mathcal S=(S_1,\cdots,S_r)$ be a nonempty ordered family of split graded-simple superalgebras over $\kk$, and let $d$ be the bound in \eqref{eq:variable-bound}. Let $\KK/\kk$ be an algebraically closed extension and $B$ a finite-dimensional $\KK$-superalgebra. Write $s=\dim_{\KK}(B/J(B))$, and suppose that $J(B)^\nu=0$. For every $L\geq ds+2\nu-1$,
the polynomial $f_{\mathcal S,L}$ from
Definition~\ref{def:separating-polynomials} does not belong to $I(\mathcal S)$.
Moreover, $\Id_{\kk}(\E_{\KK}(B))\subseteq I(\mathcal S)$
if and only if $f_{\mathcal S,L}\notin\Id_{\kk}(\E_{\KK}(B))$.
\end{theorem}
\begin{proof}
Lemma~\ref{lem:target-evaluation} and
Remark~\ref{rem:target-evaluation-field} give
$f_{\mathcal S,L}\notin I(\mathcal S)$ over $\kk$.
Hence $\Id_{\kk}(\E_{\KK}(B))\subseteq I(\mathcal S)$ implies
$f_{\mathcal S,L}\notin\Id_{\kk}(\E_{\KK}(B))$.

Conversely, suppose that $f_{\mathcal S,L}$ is not an identity of $\E_{\KK}(B)$.
Lemma~\ref{lem:sign-transfer} shows that $g_L$, the complete
linearization of $F_L$, is not a graded identity of
$B\mathbin{\widehat{\otimes}}\E_{\KK}$.
Theorem~\ref{thm:linearization} then shows that $F_L$ is not a
graded identity either.
Let $N=ds+\nu-1$ and $C=B\mathbin{\widehat{\otimes}}\E_{\KK}^{(N)}$. Since $F_L$ uses at most $d$ variables, Lemma~\ref{lem:finite-exterior} gives a nonzero graded evaluation of $F_L$ on $C$. By Lemma~\ref{lem:tensor-radical}, $
J(C)^{\nu+N}=0$ and $
\nu+N=ds+2\nu-1\leq L$.

Lemma~\ref{lem:extraction}, applied over $\KK$, therefore gives
$
\Id_{\KK}(\E_{\KK}(C))
\subseteq\prod_{i=1}^r
\Id_{\KK}(\E_{\KK}(S_i\otimes_{\kk}\KK))$.
Corollary~\ref{cor:envelope-stability} gives $\Id_{\KK}(\E_{\KK}(C))=\Id_{\KK}(\E_{\KK}(B))$.
Intersecting both sides with the free algebra over $\kk$
and applying Lemma~\ref{lem:scalar-extension} gives
$\Id_{\kk}(\E_{\KK}(B))\subseteq I(\mathcal S)$.

\end{proof}

The uniform separation criterion gives the following
intersection theorem.

\begin{theorem}\label{thm:bounded-intersection}
Let $\{T_\lambda\}_{\lambda\in\Lambda}$ be a nonempty
family of $T$-ideals. Suppose that there are integers
$s_0\geq0$ and $\nu_0\geq1$ such that every nonzero
proper $T_\lambda$ has the form
$T_\lambda=\Id_{\kk}(\E_{\KK_\lambda}(B_\lambda))$,
where $\KK_\lambda$ is an algebraically closed extension
of $\kk$ and $B_\lambda$ is a finite-dimensional
$\KK_\lambda$-superalgebra satisfying
$\dim_{\KK_\lambda}(B_\lambda/J(B_\lambda))\leq s_0$
and $J(B_\lambda)^{\nu_0}=0$.
If $P_1,\cdots,P_r$ are proper verbally prime
$T$-ideals, where $r\geq1$, then
$\bigcap_{\lambda\in\Lambda}T_\lambda
\subseteq P_1\cdots P_r$ only if
$T_\lambda\subseteq P_1\cdots P_r$
for some $\lambda\in\Lambda$.
\end{theorem}

\begin{proof}
Let $I=P_1\cdots P_r$. Assume that $T_\lambda\nsubseteq I$ for every $\lambda$. In particular, every $T_\lambda$ is nonzero.

Suppose first that every $P_i$ is nonzero.
By Theorem~\ref{thm:kemer-classification}, choose split graded-simple
superalgebras $S_i$ over $\kk$ such that
$P_i=\Id(\E(S_i))$.
Let $\mathcal S=(S_1,\cdots,S_r)$. Then $I=I(\mathcal S)$.
Let $d$ be the bound associated with $\mathcal S$ and choose $L=ds_0+2\nu_0-1$.
Theorem~\ref{thm:uniform-separation} shows that
$f_{\mathcal S,L}$ belongs to every proper $T_\lambda$.
It also belongs to $\F$, so
$
f_{\mathcal S,L}\in
(\bigcap_{\lambda\in\Lambda}T_\lambda)\setminus I$.

Suppose now that some $P_i=0$, so $I=0$.
There are only finitely many graded-isomorphism classes
of split graded-simple superalgebras of dimension at most
$s_0$. 
Choose one from each class and form the set 
$\mathcal T$. Let
$
U=\bigcap_{S\in\mathcal T}\Id_{\kk}(\E_{\kk}(S))$.
If $\mathcal T$ is empty, let $U=\F$. Since $\Id_{\kk}(\E_{\kk}(S)) \neq 0$ for each $S\in\mathcal T$ and $\F$ has no zero divisors, $U\ne0$ and $U^{\nu_0}\ne0$.

For each nonzero proper $T_\lambda$, every
graded-simple component of
$B_\lambda/J(B_\lambda)$ is graded-isomorphic
to $S\otimes_\kk \KK_\lambda$ for some
$S\in\mathcal T$, and Lemma~\ref{lem:scalar-extension} gives
$U\subseteq
\Id_\kk(\E_{\KK_\lambda}(B_\lambda/J(B_\lambda)))$. Let $
\pi_\lambda:\E_{\KK_\lambda}(B_\lambda)
\longrightarrow
\E_{\KK_\lambda}(B_\lambda/J(B_\lambda))$
be the quotient map. For every $u=u(x_1,\cdots,x_n)\in U$ and
$a_1,\cdots,a_n\in\E_{\KK_\lambda}(B_\lambda)$, $
\pi_\lambda(u(a_1,\cdots,a_n))
=u(\pi_\lambda(a_1),\cdots,\pi_\lambda(a_n))
=0$, and $u(a_1,\cdots,a_n) \in \ker\pi_\lambda=\E_{\KK_\lambda}(J(B_{\lambda}))$.
 $J(B_\lambda)^{\nu_0}=0$ gives $
(\E_{\KK_\lambda}(J(B_\lambda)))^{\nu_0}=0$.
Consequently, every polynomial in $U^{\nu_0}$
vanishes on $\E_{\KK_\lambda}(B_\lambda)$, so $
U^{\nu_0}\subseteq
\Id_\kk(\E_{\KK_\lambda}(B_\lambda))
=T_\lambda$.
This also holds for $T_\lambda=\F$, and hence
$0\ne U^{\nu_0}\subseteq\bigcap_{\lambda\in\Lambda}T_\lambda$.

\end{proof}

For finite families, Theorem~\ref{thm:bounded-intersection}
shows that Berele's Conjecture~3 in \cite{Berele2012} holds.

\begin{theorem}[\cite{Berele2012}, Conjecture~3]
\label{thm:finite-intersection}
Let $P_1,\cdots,P_r$ be proper verbally prime
$T$-ideals of $\F$, where $r\geq1$, and let
$T_1,\cdots,T_m$ be $T$-ideals, where $m\geq1$.
If $\bigcap_{j=1}^m T_j\subseteq P_1\cdots P_r$, then $T_j\subseteq P_1\cdots P_r $ for some $j$.
\end{theorem}

\begin{proof}
If some $T_j=0$, the conclusion is immediate.
By Theorem~\ref{thm:kemer-representability}, every nonzero
proper $T_j$ is the identity ideal over $\kk$ of the
Grassmann envelope of a finite-dimensional superalgebra
over an algebraically closed extension of $\kk$.
Since there are only finitely many such $T$-ideals,
the uniform bounds required in
Theorem~\ref{thm:bounded-intersection} exist.
The conclusion follows from that theorem.
\end{proof}

Berele \cite[Corollary~2.4]{Berele2012} shows that his
Conjecture~3 implies Conjecture~1. To state this consequence of
Theorem~\ref{thm:finite-intersection}, we  introduce 
the following definitions.

\begin{definition}\label{def:incidence-graph}
Let $G$ be a finite directed graph with vertex set
$\{1,\cdots,n\}$, without multiple edges but possibly
with loops. We call $G$ transitive if the existence of
edges $i\to j$ and $j\to k$ implies the existence of
an edge $i\to k$.

For a transitive graph $G$, the associated incidence
algebra is the subalgebra 
$
A(G)=\operatorname{span}_{\kk}
\{E_{ij}\mid i\to j\text{ is an edge of }G\}$ of $ \Mat_n(\kk)$,
where $E_{ij}$ denotes the standard matrix unit.

A chain in $G$ is a nonempty subset $C$ of its vertex
set such that any two distinct vertices in $C$ are
joined by an edge in at least one direction.
We identify $C$ with the subgraph of $G$ induced by $C$.

An embedding of directed graphs $G\hookrightarrow H$
is an injective map $\varphi:V(G)\hookrightarrow V(H)$ such that
$\varphi(u)\to\varphi(v)$ is an edge of $H$ whenever
$u\to v$ is an edge of $G$.
\end{definition}

The algebra $A(G)$ contains the identity matrix if
and only if every vertex of $G$ has a loop.

\begin{corollary}[\cite{Berele2012}, Conjecture~1]
\label{cor:berele-one}
Let $G_1$ and $G_2$ be finite transitive directed
graphs without multiple edges such that both
$A(G_1)$ and $A(G_2)$ contain $1$.
Then
$\Id(A(G_1))\subseteq\Id(A(G_2))$
if and only if every chain $C$ contained in $G_2$
has an embedding into $G_1$.
\end{corollary}

% ===== parts/capacity =====
\section{Inclusions between products of verbally prime
\texorpdfstring{$T$}{T}-ideals} \label{sec:capacity}

We compare two ordered products of verbally prime $T$-ideals
using the disjoint index sets constructed in
Lemma~\ref{lem:extraction}.
The criterion records the sizes and parities of these index sets
and the order of the graded-simple components. We then use it
to prove Berele's Conjecture~2 and answer his Question~2 \cite{Berele2012}.
We use $P(S)$ and $I(\mathcal S)$ as in
Definition~\ref{def:ordered-family}.
Since $\Mat_{p|q}(\kk)\cong\Mat_{q|p}(\kk)$
as graded algebras, we assume $p\geq q\geq0$
for every matrix-type component.
The ordered family may contain both matrix-type components
and components of type $\Qtype_n$, with repetitions allowed.

We first describe the index sets used in the four
constructions in the proof of Lemma~\ref{lem:extraction}.
For $\Mat_{p|q}$, an index refers to a row and the
corresponding column in its standard matrix representation.
Its parity is the parity of the corresponding basis vector.
For $\Qtype_m$, we use the $m$ indices of its even
matrix algebra $\Mat_m$, rather than the $2m$ indices
of its realization in $\Mat_{m|m}$.

\begin{lemma}\label{lem:local-capacity}
The constructions in the proof of
Lemma~\ref{lem:extraction} satisfy the following
index requirements.  Conversely, any choice of distinct indices satisfying
these requirements gives a graded embedding of the
first algebra into the second.

\begin{enumerate}
\item
For $\Mat_{a|b}$ in $\Mat_{p|q}$, the construction
uses $a$ even and $b$ odd indices of $\Mat_{p|q}$,
or $b$ even and $a$ odd indices of $\Mat_{p|q}$.
\item
For $\Mat_{a|b}$ in $\Qtype_m$, the construction
uses $a+b$ indices of the even part $\Mat_m$
of $\Qtype_m$.
\item
For $\Qtype_n$ in $\Mat_{p|q}$, the construction
uses $n$ even and $n$ odd indices of $\Mat_{p|q}$.
\item
For $\Qtype_n$ in $\Qtype_m$, the construction
uses $n$ indices of the even part $\Mat_m$
of $\Qtype_m$.
\end{enumerate}
Within a fixed algebra $\Mat_{p|q}$ or $\Qtype_m$,
graded subalgebras constructed from pairwise disjoint
index sets have pairwise orthogonal identity elements.

\end{lemma}
\begin{proof}
For the first two cases, let
$c_1,\cdots,c_{a+b}\in\Ztwo$ be the parities
of the standard basis vectors defining $\Mat_{a|b}$.

For the first case, choose distinct indices
$r_1,\cdots,r_{a+b}$ in $\Mat_{p|q}$ whose parities
are $c_1,\cdots,c_{a+b}$, or are obtained from
these by changing every parity.
The map
defined by $
E_{ij}\longmapsto E_{r_i,r_j}$
is injective and preserves the matrix-unit
multiplication rules.
Changing every parity leaves the degree
$c_i+c_j$ of each matrix unit unchanged,
so this map also preserves the grading.

In the layer $\ell_0$ chosen in the proof of
Lemma~\ref{lem:extraction}, the row indices recorded
for the variables $x_{ij}$ are pairwise distinct
because $H(\lambda^{(\ell_0)})\ne0$.
For $1\leq j<a_i$, the column index of the matrix
unit chosen for $x_{ij}$ agrees with the row index
of the matrix unit chosen for $x_{i,j+1}$,
since their product is nonzero.
The given degrees
$|x_{ij}|=c_{ij}+c_{i,j+1}$ therefore determine
the parities of the recorded row indices up to
one common change of parity.
Thus these indices satisfy the stated requirement.

For the second case, choose distinct indices
$r_1,\cdots,r_{a+b}$ of the even matrix algebra
of $\Qtype_m$. Define
 a map by $
E_{ij}\longmapsto u^{c_i+c_j}E_{r_i,r_j}$.
The images are linearly independent and have
degrees $c_i+c_j$.
Since $u^2=1$ and $u$ commutes with every
matrix unit, they satisfy the matrix-unit
multiplication rules. Hence this map is
a graded embedding.

In the proof of Lemma~\ref{lem:extraction},
the row indices recorded for the variables $x_{ij}$
are indices of the even matrix algebra $\Mat_m$
of $\Qtype_m$.
They are pairwise distinct because
$H(\lambda^{(\ell_0)})\ne0$.

For the third case, choose distinct even
indices $r_1,\cdots,r_n$ and distinct odd indices
$s_1,\cdots,s_n$ in $\Mat_{p|q}$.
Define
$
A_{ij}=E_{r_i,r_j}+E_{s_i,s_j}$ and $
B_{ij}=E_{r_i,s_j}+E_{s_i,r_j}$.
The elements $A_{ij}$ are even, and the elements
$B_{ij}$ are odd. They satisfy $
A_{ij}A_{kl}=B_{ij}B_{kl}=\delta_{jk}A_{il}$ and $
A_{ij}B_{kl}=B_{ij}A_{kl}=\delta_{jk}B_{il}$.
Their linear independence therefore shows that the map defined by $
E_{ij}\longmapsto A_{ij}$ and $
uE_{ij}\longmapsto B_{ij}$
is a graded embedding of $\Qtype_n$
into $\Mat_{p|q}$.

In the layer $\ell_0$ chosen in the proof of
Lemma~\ref{lem:extraction}, the matrix units chosen
for the variables $x^+_{ij}$ are even and have
nonzero product. Their row indices therefore
have the same parity.
The same holds for the variables $x^-_{ij}$.
The odd matrix unit chosen for $\xi_i$ joins
the two families, so their row indices have
opposite parities.
Within each set, the recorded row indices are
distinct because $H(\lambda^{(\ell_0)})\ne0$.

For the fourth case, choose distinct indices
$r_1,\cdots,r_n$ of the even matrix algebra
of $\Qtype_m$. The map defined by $
E_{ij}\longmapsto E_{r_i,r_j}$ and $
uE_{ij}\longmapsto uE_{r_i,r_j}$
is a graded embedding of $\Qtype_n$
into $\Qtype_m$.

In the proof of Lemma~\ref{lem:extraction},
the row indices recorded for the variables $x^+_{ij}$
are pairwise distinct because
$H(\lambda^{(\ell_0)})\ne0$.
These indices therefore give the required embedding.
The indices recorded for the variables $x^-_{ij}$
may overlap with this set but are not needed
for the construction.

In each case, the unit of the image is the sum
of the diagonal matrix units indexed by the
chosen set.
For the embedding of $\Qtype_n$ into $\Mat_{p|q}$,
this unit is $
\sum_{i=1}^n
(E_{r_i,r_i}+E_{s_i,s_i})$.
Thus the units of subalgebras constructed from
disjoint index sets multiply to zero in either
order.
\end{proof}

To compare ordered products of verbally prime $T$-ideals,
we combine the index requirements in
Lemma~\ref{lem:local-capacity} with the order of the
graded-simple components.
Several consecutive members of one ordered family
may be embedded into the same component of another
ordered family using pairwise disjoint index sets.

\begin{definition}\label{def:super-capacity}
Let $\mathcal S=(S_1,\cdots,S_r)$ and $\mathcal R=(R_1,\cdots,R_s)$ be nonempty ordered families of split graded-simple superalgebras. Write $\mathcal S\preceq_{\mathrm{sup}}\mathcal R$ if there is a nondecreasing map $\phi:\{1,\cdots,r\}\to\{1,\cdots,s\}$ with the following properties:

\begin{enumerate}
\item
If $R_j=\Mat_{p_j|q_j}$, choose, independently for each $S_i=\Mat_{a_i|b_i}$ with $\phi(i)=j$, a pair $(\alpha_i,\beta_i)$ equal to either $(a_i,b_i)$ or $(b_i,a_i)$. These choices must satisfy
\[
\sum_{\substack{\phi(i)=j\\ S_i=\Mat_{a_i|b_i}}}\alpha_i
+\sum_{\substack{\phi(i)=j\\ S_i=\Qtype_{n_i}}}n_i\leq p_j,
\qquad
\sum_{\substack{\phi(i)=j\\ S_i=\Mat_{a_i|b_i}}}\beta_i
+\sum_{\substack{\phi(i)=j\\ S_i=\Qtype_{n_i}}}n_i\leq q_j.
\]
\item
If $R_j=\Qtype_{m_j}$, the requirement is
\[
\sum_{\substack{\phi(i)=j\\ S_i=\Mat_{a_i|b_i}}}(a_i+b_i)
+\sum_{\substack{\phi(i)=j\\ S_i=\Qtype_{n_i}}}n_i\leq m_j.
\]

\end{enumerate}
The map $\phi$ is not required to be surjective.
If no $i$ satisfies $\phi(i)=j$, the sums in the
inequalities for $R_j$ are empty and hence zero,
so these inequalities hold.

\end{definition}

The relation just defined characterizes inclusions between
the corresponding products of verbally prime $T$-ideals.

\begin{theorem}\label{thm:mixed-capacity}
Let $\mathcal S$ and $\mathcal R$ be nonempty ordered families of split graded-simple superalgebras over a field of characteristic zero. Then
$
I(\mathcal R)\subseteq I(\mathcal S)$ if and only if
 $\mathcal S\preceq_{\mathrm{sup}}\mathcal R$.

\end{theorem}
\begin{proof}

We first work over an algebraically closed field.
Let $B=A(\mathcal R)$ be the superalgebra in
Lemma~\ref{lem:target-evaluation}; it is minimal.
For each diagonal block $R_j$, let $d_j=p_j+q_j$ if
$R_j=\Mat_{p_j|q_j}$, and let $d_j=2n_j$ if $R_j=\Qtype_{n_j}$.
The strictly upper $(i,j)$-block of $B$ is
$\Mat_{d_i\times d_j}(\kk)$ for $i<j$, and
$\Id(\E(B))=I(\mathcal R)$.

Suppose that $I(\mathcal R)\subseteq I(\mathcal S)$. Let $d$ be the variable bound for the polynomials associated with $\mathcal S$, let $s_0=\dim(B/J(B))$, and choose $\nu$ with $J(B)^\nu=0$. Choose $N\geq ds_0+\nu-1$ and let $C=B\mathbin{\widehat{\otimes}}\E^{(N)}$. Choose $L$ such that $J(C)^L=0$.  The polynomial $f_{\mathcal S,L}$ does not belong to $I(\mathcal S)$, so it does not belong to $\Id(\E(B))$. By Lemma~\ref{lem:sign-transfer} and Theorem~\ref{thm:linearization}, $F_L$ is not a graded identity of $B\mathbin{\widehat{\otimes}}\E$. Lemma~\ref{lem:finite-exterior} then gives a nonzero graded evaluation on $C$.

The algebra $C$ has
semisimple part $(R_1\oplus\cdots\oplus R_s)\otimes1$.
Apply Lemma~\ref{lem:extraction} with this  part.
For each $i$, it constructs $S'_i\cong S_i$ inside one
$R_j\otimes1$; let $\phi(i)=j$.
Within each $R_j\otimes1$, the index sets used for distinct
$S'_i$ are disjoint.
Lemma~\ref{lem:local-capacity} therefore gives the inequalities
in Definition~\ref{def:super-capacity}.

It remains to check that $\phi$ is nondecreasing.
 $B$ is block upper triangular with diagonal blocks
$R_1,\cdots,R_s$, and
$C$ has zero $(a,b)$-blocks whenever $a>b$.
By Lemma~\ref{lem:extraction}, there are homogeneous idempotents
$e_i\in S'_i$ such that $e_iCe_{i+1}\ne0$ for $1\leq i<r$.
Since $e_i\in R_{\phi(i)}\otimes1$,  $e_iCe_{i+1}$
is contained in the $(\phi(i),\phi(i+1))$-block of $C$.
Hence $\phi(i)\leq\phi(i+1)$ for every $1\leq i<r$,
so $\phi$ is nondecreasing.

Conversely, suppose that $\mathcal S\preceq_{\mathrm{sup}}\mathcal R$.
In each $R_j$, choose disjoint index sets satisfying the
inequalities in Definition~\ref{def:super-capacity}.
Lemma~\ref{lem:local-capacity} gives
$S'_i\cong S_i$ with pairwise orthogonal units $\epsilon_i$.
For each $i$, choose a minimal nonzero homogeneous diagonal
idempotent $e_i\in S'_i$ and the index of one of its
nonzero diagonal entries.

Let $\rho_i$ denote the selected index for $e_i$ in the
matrix representation of $B$, and let $\widetilde E_{ab}$
denote the matrix units of the full matrix algebra
containing $B$.
If $S'_i\subseteq R_j=\Qtype_m$, choose $\rho_i$ among
the first $m$ row indices of this diagonal block.
This is possible because the two diagonal $m\times m$
blocks of $e_i$ are equal and nonzero.

We construct homogeneous elements $v_i\in e_iBe_{i+1}$.
When $S'_i$ and $S'_{i+1}$ lie in the same component
$R_j=\Mat_{p_j|q_j}$, set
$v_i=\widetilde E_{\rho_i,\rho_{i+1}}$.
When they lie in the same component $R_j=\Qtype_m$, set
$v_i=\widetilde E_{\rho_i,\rho_{i+1}}
+\widetilde E_{\rho_i+m,\rho_{i+1}+m}$,
which is an even element of $R_j$.
When they lie in distinct components, set
$v_i=\widetilde E_{\rho_i,\rho_{i+1}}$.
Since $\phi(i)<\phi(i+1)$, this matrix unit belongs
to the full $(\phi(i),\phi(i+1))$-block of $B$.
In every case, $v_i$ is homogeneous and belongs to
$e_iBe_{i+1}$.
The $(\rho_1,\rho_r)$-entry of $v_1\cdots v_{r-1}$
is one, so $v_1\cdots v_{r-1}\ne0$.

Let $C'$ be the graded subalgebra of $B$ generated by the $S'_i$ and the $v_i$, and let $R'$ be its ideal generated by the $v_i$.
Since $S'_i\subseteq\epsilon_iB\epsilon_i$ and
$v_i\in\epsilon_iB\epsilon_{i+1}$, orthogonality of
the $\epsilon_i$ implies that every nonzero product
of these generators containing at least one $v_i$
lies in $\epsilon_aB\epsilon_b$ for some $a<b$.
 Consequently,
$
R'\subseteq\bigoplus_{i<j}\epsilon_iB\epsilon_j$, $(R')^r=0$, and
$C'=(S'_1\oplus\cdots\oplus S'_r)\oplus R'$.
Since $R'$ is nilpotent and $C'/R'$ is semisimple, $R'=J(C')$. Moreover, $v_i=e_iv_ie_{i+1}$ and their product is nonzero. Thus $C'$ is minimal. For $r=1$, take $C'=S'_1$.  By Theorem~\ref{thm:minimal-product}, we have $
I(\mathcal R)=\Id(\E(B))
\subseteq\Id(\E(C'))=I(\mathcal S)$.

For an arbitrary field $\kk$ of characteristic zero, extend
scalars to an algebraic closure.
The inequalities in Definition~\ref{def:super-capacity} are unchanged,
and Lemma~\ref{lem:scalar-extension} preserves the identity ideals
and their products and recovers them by intersection with $\F$.
So, the equivalence  holds over $\kk$.
\end{proof}

We now apply Theorem~\ref{thm:mixed-capacity} to
Conjecture~2 and Question~2 of \cite{Berele2012}.
To compare the notation, observe that
$P(\Mat_{n|0})=\Id(\Mat_n(\kk))$ and
$P(\Qtype_n)=\Id(\Mat_n(\E))$.
Indeed, $\E(\Mat_{n|0})=\Mat_n(\E_0)$ has the same
multilinear identities as $\Mat_n(\kk)$, since $\E_0$
is commutative and the constant matrices form a subalgebra.
Also, $\E(\Qtype_n)\cong\Mat_n(\E)$ as ordinary algebras.
For $p,q\geq1$, the envelope $\E(\Mat_{p|q})$ is the
algebra denoted by $M_{p,q}$ in \cite{Berele2012}.

We first determine all inclusions between proper
verbally prime $T$-ideals, thereby proving Berele’s Conjecture 2.

\begin{corollary}[\cite{Berele2012}, Conjecture~2]
\label{cor:prime-inclusions}
Let $U$ and $V$ be proper verbally prime $T$-ideals
of $\F$.
Then $U\subseteq V$ if and only if $U=0$, or
$U=P(R)$ and $V=P(S)$ for split graded-simple
superalgebras $R$ and $S$ over $\kk$ satisfying one of the
following conditions:
\begin{enumerate}
\item
$R=\Mat_{p|q}$ and $S=\Mat_{a|b}$, where
$(a\leq p,\ b\leq q)$ or $(a\leq q,\ b\leq p)$.
\item
$R=\Mat_{p|q}$ and $S=\Qtype_n$, where
$n\leq\min\{p,q\}$.
\item
$R=\Qtype_m$ and $S=\Mat_{a|b}$, where
$a+b\leq m$.
\item
$R=\Qtype_m$ and $S=\Qtype_n$, where $n\leq m$.
\end{enumerate}
Here $p,q,a,b\geq0$, with $p+q>0$ and $a+b>0$,
and $m,n\geq1$.
The seven inclusions listed before Conjecture~2
in \cite{Berele2012} generate all inclusions between
nonzero proper verbally prime $T$-ideals.
\end{corollary}

\begin{proof}
It suffices to consider the case $U,V\ne0$.
By Theorem~\ref{thm:kemer-classification} and
Lemma~\ref{lem:scalar-extension}, we can write
$U=P(R)$ and $V=P(S)$ for split graded-simple
superalgebras $R$ and $S$ over $\kk$.
Each algebra is of matrix type or of type $\Qtype$,
so the four cases above exhaust all possibilities.
Applying Theorem~\ref{thm:mixed-capacity} with
$\mathcal R=(R)$ and $\mathcal S=(S)$ gives
exactly the inequalities in \textup{(i)--(iv)}.

It remains to show that every inclusion allowed by
these inequalities follows from Berele's seven
inclusions.
Case \textup{(i)} follows from Berele's first, third,
and fourth inclusions, and case \textup{(ii)} from
his fourth and seventh.
Case \textup{(iii)} follows from his first, second,
fifth, and sixth inclusions, and case \textup{(iv)}
from his sixth.
\end{proof}

Question~2 of \cite{Berele2012} asks when a proper
verbally prime $T$-ideal is contained in a product
of such ideals, and singles out the case of one
factor on the right.
Corollary~\ref{cor:prime-inclusions} settles this
special case.
We now answer the general question by taking
$\mathcal R$ to have one member in
Theorem~\ref{thm:mixed-capacity}.

\begin{corollary}[Answer to Question~2 of \cite{Berele2012}]
\label{cor:question-two}

Let $\mathcal S=(S_1,\cdots,S_r)$ be a nonempty ordered
family of split graded-simple superalgebras.
Let $p,q\geq0$ with $p+q>0$, and let $m\geq1$.
\begin{enumerate}
\item
The inclusion
$P(\Mat_{p|q})\subseteq I(\mathcal S)$
holds if and only if, for each $i$ with
$S_i=\Mat_{a_i|b_i}$, one can choose
$(\alpha_i,\beta_i)=(a_i,b_i)$ or $(b_i,a_i)$
so that
\[
\sum_{S_i=\Mat_{a_i|b_i}}\alpha_i
+\sum_{S_i=\Qtype_{n_i}}n_i\leq p,
\qquad
\sum_{S_i=\Mat_{a_i|b_i}}\beta_i
+\sum_{S_i=\Qtype_{n_i}}n_i\leq q.
\]

\item
The inclusion
$P(\Qtype_m)\subseteq I(\mathcal S)$
holds if and only if
\[
\sum_{S_i=\Mat_{a_i|b_i}}(a_i+b_i)
+\sum_{S_i=\Qtype_{n_i}}n_i\leq m.
\]
\end{enumerate}
\end{corollary}

\begin{proof}
Apply Theorem~\ref{thm:mixed-capacity} with
$\mathcal R=(\Mat_{p|q})$ or $\mathcal R=(\Qtype_m)$,
respectively.
The map $\phi$ in Definition~\ref{def:super-capacity}
must send every index to $1$, so the inequalities
in that definition become the displayed conditions.
\end{proof}

The conditions in Corollary~\ref{cor:question-two}
are unchanged by permuting $S_1,\cdots,S_r$,
since all these algebras are assigned to the same
component in Definition~\ref{def:super-capacity}.

The following example illustrates
Corollary~\ref{cor:question-two} and shows that the
dimensions of the graded-simple components alone
do not determine the inclusion.

\begin{example}\label{ex:mixed-capacity}
Let $\mathcal S=(\Mat_{2|1},\Qtype_1)$.
In $\Mat_{3|2}$, choose two even and one odd index for
$\Mat_{2|1}$, and the remaining even and odd indices
for $\Qtype_1$.
Corollary~\ref{cor:question-two} gives
$P(\Mat_{3|2})\subseteq P(\Mat_{2|1})P(\Qtype_1)$.

In $\Mat_{4|1}$, the construction of $\Qtype_1$ uses
the only odd index, whereas either choice for
$\Mat_{2|1}$ requires at least one further odd index.
Hence
$P(\Mat_{4|1})\not\subseteq P(\Mat_{2|1})P(\Qtype_1)$,
although $\Mat_{3|2}$ and $\Mat_{4|1}$ both have
dimension $25$.

Inside $\Qtype_m$, the constructions of $\Mat_{2|1}$ and
$\Qtype_1$ require three indices and one index, respectively,
of the even subalgebra $\Mat_m$.
Corollary~\ref{cor:question-two} therefore shows that
$
P(\Qtype_m)\subseteq P(\Mat_{2|1})P(\Qtype_1)
$ if and only if $m\geq4$.

\end{example}

\section{Applications to intersections of products
of verbally prime \texorpdfstring{$T$}{T}-ideals}
\label{sec:antichains}

We apply Theorem~\ref{thm:bounded-intersection}
and Theorem~\ref{thm:mixed-capacity} to intersections of products
of nonzero proper verbally prime $T$-ideals.
The sets indexing these intersections need not be finite.
The degree of a nonzero polynomial in the intersection
bounds both the cardinalities of the ordered families
and the parameters of their members.
Consequently, only finitely many distinct ordered families
can occur. These bounds also allow us to apply
Theorem~\ref{thm:bounded-intersection} to characterize
inclusions between such intersections.
We then determine which products can be omitted from
a nonzero intersection and prove that the set of
remaining products is uniquely determined by the
intersection. 
The final example shows why the case in which the intersection is zero must be treated separately in the
inclusion criterion.
Throughout this section, each member of a nonempty finite ordered family
is either $\Mat_{a|b}(\kk)$ with $a\geq b\geq0$ and $a+b>0$,
or $\Qtype_n$ with $n\geq1$.
We use $P(S)$ and $I(\mathcal S)$ as in
Definition~\ref{def:ordered-family}.

\begin{definition}\label{def:family-intersection}
For a nonempty set $\Gamma$ of nonempty finite  ordered families, let
$\mathcal I_\Gamma=\bigcap_{\gamma\in\Gamma}I(\gamma)$.
The factors in $I(\gamma)$ occur in the order of the members
of $\gamma$.
The set $\Gamma$ may be infinite, although each of its
members is a finite ordered family.
\end{definition}

We first bound the ordered families that can occur when
the intersection contains a nonzero polynomial.

\begin{lemma}\label{lem:intersection-degree-bounds}

Let $\Gamma$ be a nonempty set of nonempty finite ordered families.
Then $\mathcal I_\Gamma\ne0$ if and only if $\Gamma$ is finite.
More precisely, suppose that $0\ne f\in\mathcal I_\Gamma$,
let $D=\deg f$,
and let $q=\lfloor D/2\rfloor$.
For every $\gamma=(S_1,\cdots,S_r)\in\Gamma$, we have
$r\leq q$.
If $S_i=\Mat_{a|b}$, then $a,b\leq q$;
if $S_i=\Qtype_n$, then $n\leq q$.
\end{lemma}

\begin{proof}
Suppose that $\mathcal I_\Gamma\ne0$, and choose $f$, $D$,
and $q$ as in the statement.
The algebra $\E(S)$ is nonzero, so it has no nonzero
polynomial identity of degree one. Since every homogeneous component of
a polynomial identity is again an identity,
 $P(S)\subseteq\F^2$ for every member $S$ under
consideration.
For $\gamma=(S_1,\cdots,S_r)\in\Gamma$, it follows that
$I(\gamma)\subseteq\F^{2r}$.
Since $f\in I(\gamma)$ and $\deg f=D$, we obtain
$2r\leq D$, and therefore $r\leq q$.

 $I(\gamma)$ is contained in each 
 $P(S_i)$, so $f$ is an identity of every
$\E(S_i)$.
If $S_i=\Mat_{a|b}$, then $\E(S_i)$ contains a subalgebra
isomorphic to $\Mat_a(\kk)$ when $a>0$, and a subalgebra
isomorphic to $\Mat_b(\kk)$ when $b>0$. If $S_i=\Qtype_n$, then $\E(S_i)$ contains a subalgebra
isomorphic to $\Mat_n(\kk)$. The degree bound in Theorem~\ref{thm:amitsur-levitzki}
 gives $a,b\leq q$ and $n\leq q$, respectively.
There are only finitely many choices of these parameters,
and the cardinality of every ordered family is at most $q$.
Thus $\Gamma$ is finite.

Conversely, suppose that $\Gamma$ is finite.
Every $P(S)$ is nonzero, and $\F$ has no zero
divisors, so every product $I(\gamma)$ is nonzero.
The product of the finitely many ideals $I(\gamma)$ is nonzero and is contained
in each of them.
Hence $\mathcal I_\Gamma\ne0$.
\end{proof}

When the intersection is nonzero, the bounds in
Lemma~\ref{lem:intersection-degree-bounds} allow us
to apply Theorem~\ref{thm:bounded-intersection}.

\begin{theorem}\label{thm:intersection-capacity}
For nonempty sets $\Gamma$ and $\Delta$ of nonempty finite ordered
families, $\mathcal I_\Gamma\subseteq\mathcal I_\Delta$
if and only if $\Gamma$ is infinite or, for every
$\delta\in\Delta$, there is $\gamma\in\Gamma$ with
$\delta\preceq_{\mathrm{sup}}\gamma$.
\end{theorem}

\begin{proof}
If $\Gamma$ is infinite, then
$\mathcal I_\Gamma=0$ by
Lemma~\ref{lem:intersection-degree-bounds}, so the
inclusion holds.

Suppose that $\mathcal I_\Gamma\ne0$ and
$\mathcal I_\Gamma\subseteq\mathcal I_\Delta$.
Choose $0\ne f\in\mathcal I_\Gamma$ and let
$q=\lfloor\deg f/2\rfloor$.
Lemma~\ref{lem:intersection-degree-bounds} gives $q\geq1$
and bounds the cardinality and all parameters of every
$\gamma\in\Gamma$ by $q$.

Let $\KK$ be an algebraic closure of $\kk$.
For $\gamma=(S_1,\cdots,S_r)\in\Gamma$, let $B_\gamma$
be the block upper triangular superalgebra of
Lemma~\ref{lem:target-evaluation}, constructed over $\KK$
with diagonal blocks $S_i\otimes_\kk \KK$.
By Lemma~\ref{lem:target-evaluation}
and Lemma~\ref{lem:scalar-extension},
$I(\gamma)=\Id_\kk(\E_{\KK}(B_\gamma))$.
 $J(B_\gamma)$ consists of the strictly upper blocks,
so $J(B_\gamma)^r=0$ and hence $J(B_\gamma)^q=0$.
Each diagonal block has dimension at most $4q^2$,
and there are at most $q$ blocks.
Thus $\dim_{\KK}(B_\gamma/J(B_\gamma))\leq4q^3$.

Let $\delta\in\Delta$.
We have $\bigcap_{\gamma\in\Gamma}I(\gamma)
\subseteq I(\delta)$.
Theorem~\ref{thm:bounded-intersection}, with
$s_0=4q^3$ and $\nu_0=q$, gives $\gamma\in\Gamma$
such that $I(\gamma)\subseteq I(\delta)$.
Theorem~\ref{thm:mixed-capacity} then gives
$\delta\preceq_{\mathrm{sup}}\gamma$.

Conversely, suppose that for every $\delta\in\Delta$
there is $\gamma\in\Gamma$ with
$\delta\preceq_{\mathrm{sup}}\gamma$.
Theorem~\ref{thm:mixed-capacity} gives
$\mathcal I_\Gamma\subseteq I(\gamma)\subseteq I(\delta)$.
Intersecting over $\delta\in\Delta$ proves
$\mathcal I_\Gamma\subseteq\mathcal I_\Delta$.
\end{proof}

To remove redundant products from a nonzero intersection,
we use the partial order established below.

\begin{lemma}\label{lem:capacity-order}

The relation $\preceq_{\mathrm{sup}}$ is a partial order
on nonempty finite ordered families.
% whose members are
%$\Mat_{a|b}$ with $a\geq b\geq0$ and $a+b>0$,
%or $\Qtype_n$ with $n\geq1$.

Let
$
\rho(\Mat_{a|b})=(a+b)^2$,
$\rho(\Qtype_n)=2n^2$, and $
\rho(\mathcal S)=\sum_i\rho(S_i)$.
If $\mathcal S\preceq_{\mathrm{sup}}\mathcal R$, then
$\rho(\mathcal S)\leq\rho(\mathcal R)$, with equality
if and only if $\mathcal S=\mathcal R$.

\end{lemma}
\begin{proof}
Suppose that $\mathcal S\preceq_{\mathrm{sup}}\mathcal R$,
and let $\phi$ be a nondecreasing map satisfying
Definition~\ref{def:super-capacity}.
By Definition~\ref{def:super-capacity} and
Lemma~\ref{lem:local-capacity}, choosing pairwise disjoint
index sets for the $S_i$ assigned to each $R_j$
gives a graded embedding
$
\bigoplus_{i=1}^r S_i
\hookrightarrow
\bigoplus_{j=1}^s R_j$.
Since $\rho$ is the dimension of the corresponding
direct sum, we obtain
$\rho(\mathcal S)\leq\rho(\mathcal R)$.

If equality holds, the embedding is surjective.
For each $j$, the images of the $S_i$ with
$\phi(i)=j$ then form a direct sum decomposition
of $R_j$ into graded ideals, since distinct images
have orthogonal identity elements.
As $R_j$ is graded-simple, exactly one such image
occurs, and it equals $R_j$.
Consequently, $\phi$ is a bijection.
Since $\phi$ is nondecreasing, it preserves the order,
so $r=s$ and $S_i\cong R_i$ for every $i$.
Since matrix-type components are written with
$a\geq b$ for $\Mat_{a|b}$, these graded isomorphisms
imply $\mathcal S=\mathcal R$.

Reflexivity is immediate.
If $\mathcal S\preceq_{\mathrm{sup}}\mathcal R$
and $\mathcal R\preceq_{\mathrm{sup}}\mathcal T$,
then Theorem~\ref{thm:mixed-capacity} gives
$
I(\mathcal T)\subseteq I(\mathcal R)
\subseteq I(\mathcal S)$,
and hence $\mathcal S\preceq_{\mathrm{sup}}\mathcal T$.
Thus the relation is transitive.
Finally, if $\mathcal S$ and $\mathcal R$ are related
in both directions, the dimension inequalities give
$\rho(\mathcal S)=\rho(\mathcal R)$, and therefore
$\mathcal S=\mathcal R$.
This proves antisymmetry.
\end{proof}

By Theorem~\ref{thm:mixed-capacity} and
Lemma~\ref{lem:capacity-order}, two nonempty finite ordered families
$\gamma$ and $\eta$ satisfy $I(\gamma)=I(\eta)$
if and only if $\gamma=\eta$.

\begin{definition}\label{def:reduced-antichain}
For a finite nonempty set $\Gamma$ of such ordered
families, let $\Gamma_{\max}$ be the set of its maximal
members with respect to $\preceq_{\mathrm{sup}}$.
\end{definition}

\begin{theorem}\label{thm:unique-antichain}
Let $\Gamma$ and $\Delta$ be nonempty sets of nonempty finite
ordered families.
If $\Gamma$ is finite, then
$\mathcal I_\Gamma=\mathcal I_{\Gamma_{\max}}$.
Moreover, $\mathcal I_\Gamma=\mathcal I_\Delta$ if and
only if either both sets are infinite, or both are
finite and $\Gamma_{\max}=\Delta_{\max}$.

Let $\Lambda$ be a nonempty set.
For each $\lambda\in\Lambda$, let $r_\lambda\geq1$
be an integer, and let
$P_{\lambda,1},\cdots,P_{\lambda,r_\lambda}$
be nonzero proper verbally prime $T$-ideals.
Let
$
T_\lambda=P_{\lambda,1}\cdots P_{\lambda,r_\lambda}$ and $
\mathfrak a=\bigcap_{\lambda\in\Lambda}T_\lambda$.
If $\mathfrak a\ne0$, then there is a finite nonempty
subset $\Lambda_0\subseteq\Lambda$ such that
$\mathfrak a=\bigcap_{\lambda\in\Lambda_0}T_\lambda$
and
$\mathfrak a\subsetneq
\bigcap_{\mu\in\Lambda_0\setminus\{\lambda\}}T_\mu$
for every $\lambda\in\Lambda_0$.
The set $\{T_\lambda\mid\lambda\in\Lambda_0\}$
is uniquely determined by $\mathfrak a$.
If $\mathfrak a=0$, then
$
\bigcap_{\lambda\in\Lambda\setminus\{\lambda_0\}}
T_\lambda=0$
for every $\lambda_0\in\Lambda$.
Here the empty intersection is defined to be $\F$.

\end{theorem}

\begin{proof}
Suppose first that $\Gamma$ is finite.
If $\gamma\preceq_{\mathrm{sup}}\eta$, then
$I(\eta)\subseteq I(\gamma)$ by
Theorem~\ref{thm:mixed-capacity}.
Every $\gamma\in\Gamma$ satisfies
$\gamma\preceq_{\mathrm{sup}}\eta$ for some
$\eta\in\Gamma_{\max}$, so
$\mathcal I_\Gamma=\mathcal I_{\Gamma_{\max}}$.

By Lemma~\ref{lem:intersection-degree-bounds},
$\mathcal I_\Gamma=0$ if and only if $\Gamma$ is
infinite, and the same holds for $\Delta$.
Thus equality holds when both sets are infinite and
fails when precisely one is infinite.

It remains to consider finite $\Gamma$ and $\Delta$.
Suppose that $\mathcal I_\Gamma=\mathcal I_\Delta$,
and replace both sets by their sets of maximal members.
For $\gamma\in\Gamma$,
Theorem~\ref{thm:intersection-capacity} gives
$\delta\in\Delta$ with
$\gamma\preceq_{\mathrm{sup}}\delta$.
Applying the theorem in the opposite direction gives
$\gamma'\in\Gamma$ with
$\delta\preceq_{\mathrm{sup}}\gamma'$.
Thus $\gamma\preceq_{\mathrm{sup}}\gamma'$,
and the maximality of $\gamma$ implies
$\gamma=\gamma'$.
Antisymmetry  gives $\gamma=\delta$.
Hence $\Gamma\subseteq\Delta$. Similarly, $\Delta \subseteq \Gamma$.
The converse is immediate.

For finite $\Gamma$, we next show that
$\mathcal I_\Gamma\subsetneq
\bigcap_{\eta\in\Gamma_{\max}\setminus\{\gamma_0\}}I(\eta)$
for every $\gamma_0\in\Gamma_{\max}$.
If $\Gamma_{\max}=\{\gamma_0\}$, this follows from
$I(\gamma_0)\subsetneq\F$.
Otherwise, suppose that the two intersections are equal.
Theorem~\ref{thm:intersection-capacity}, applied to
$\Gamma_{\max}\setminus\{\gamma_0\}$ and
$\{\gamma_0\}$, gives
$\eta\in\Gamma_{\max}\setminus\{\gamma_0\}$ with
$\gamma_0\preceq_{\mathrm{sup}}\eta$,
contrary to the maximality of $\gamma_0$.
Since any nonmaximal member can be omitted,
the preceding equality criterion shows that
$\{I(\gamma)\mid\gamma\in\Gamma_{\max}\}$
is uniquely determined by $\mathcal I_\Gamma$.

By Theorem~\ref{thm:kemer-classification}, choose
$S_{\lambda,i}$ among the superalgebras considered
in this section such that
$P_{\lambda,i}=P(S_{\lambda,i})$.
Let
$\gamma_\lambda=(S_{\lambda,1},\cdots,S_{\lambda,r_\lambda})$
and $\Gamma=\{\gamma_\lambda\mid\lambda\in\Lambda\}$.
Then $T_\lambda=I(\gamma_\lambda)$ and
$\mathfrak a=\mathcal I_\Gamma$.
If $\mathfrak a\ne0$, then $\Gamma$ is finite by
Lemma~\ref{lem:intersection-degree-bounds}.
For each $\gamma\in\Gamma_{\max}$, choose one index
$\lambda\in\Lambda$ with $\gamma_\lambda=\gamma$,
and let $\Lambda_0$ be the set of these indices.
Then $\Lambda_0$ is finite and nonempty.
The preceding argument gives the required equality
and strict inclusions, and shows that
$\{T_\lambda\mid\lambda\in\Lambda_0\}$
is uniquely determined by $\mathfrak a$.

If $\mathfrak a=0$, then $\Gamma$ is infinite by
Lemma~\ref{lem:intersection-degree-bounds}.
Choose $\lambda_0\in\Lambda$.
Removing this index removes at most one member of
$\Gamma$, so
$\{\gamma_\lambda\mid
\lambda\in\Lambda\setminus\{\lambda_0\}\}$
is still infinite.
Applying the same lemma gives
$\bigcap_{\lambda\in\Lambda\setminus\{\lambda_0\}}
T_\lambda=0$.

\end{proof}

The following example shows why the case of a zero intersection in
Theorem~\ref{thm:intersection-capacity} must be treated
separately.

\begin{example}\label{ex:infinite-intersection-zero}
For $n\geq1$, let $\gamma_n$ be the ordered family
with $n$ members, all equal to $\Mat_{1|0}$, and let
$\Gamma=\{\gamma_n\mid n\geq1\}$.
Then $I(\gamma_n)=I_1^n$ and
$\mathcal I_\Gamma=0\subseteq I_2$, but
$I_1^n\nsubseteq I_2$ for every $n\geq1$.
Moreover, $\Gamma$ has no maximal member with respect
to $\preceq_{\mathrm{sup}}$.
\end{example}

\begin{proof}
The equality $\mathcal I_\Gamma=0$ follows from
Lemma~\ref{lem:intersection-degree-bounds}.
The polynomial
$f_n=[x_1,x_2]\cdots[x_{2n-1},x_{2n}]$
belongs to $I_1^n$.
In $\Mat_2(\kk)$, substitute $E_{12}$ for each
$x_{2i-1}$ and $E_{21}$ for each $x_{2i}$.
The value is $(E_{11}-E_{22})^n\ne0$, so
$f_n\notin I_2$.
Moreover, it is easy to check that $\Gamma$ has no maximal member.
\end{proof}

\section{Applications to path algebras}\label{sec:path-applications}

We apply Theorem~\ref{thm:finite-intersection} and Theorem~\ref{thm:mixed-capacity}
 to 
compare the polynomial
identities of path algebras and positive powers of their
arrow ideals.
 For PI quivers, finite matrix models
give chain formulas for the path algebra and its arrow ideal.
We characterize equality between the identities of a path
algebra and those of its path ideals, and describe the
stabilization of the latter as the minimum path length
increases. We conclude with numerical criteria and examples.

% ===== parts/berele =====

% ===== parts/incidence =====

% ===== parts/paths =====
\subsection{Finite matrix models and chain identities}
\label{sec:positive-paths}

Throughout this subsection, let $Q$ be a finite PI quiver.
We construct finite-dimensional matrix models with the same
polynomial identities as the path algebra and its path ideals.
Their block
decompositions give formulas for the identity ideals of the path
algebra and its arrow ideal in terms of chains of SCCs.

We first represent paths by matrix units with the same endpoints.

\begin{definition}\label{def:path-matrix-models}
Choose a numbering of $Q_0$, and let $E_{uv}$ denote the
corresponding matrix units in $\Mat_{|Q_0|}(\kk)$.
For $m\ge1$, define
\begin{align*}
 A(Q)&=\operatorname{span}_{\kk}
 \{E_{uv}\mid\text{there is a path from $u$ to $v$}\},\\
 A_m(Q)&=\operatorname{span}_{\kk}
 \{E_{uv}\mid\text{a path from $u$ to $v$ has length at least $m$}\}.
\end{align*}
The paths defining $A(Q)$ may be trivial.
Each space is
a subalgebra of $\Mat_{|Q_0|}(\kk)$.
\end{definition}

\begin{proposition}\label{prop:tail-models}
For every finite PI quiver $Q$ and every $m\ge1$,
 $\Id(\kk Q)=\Id(A(Q))$ and $
 \Id(\kk Q_{\ge m})=\Id(A_m(Q))$.
Moreover, $A_m(Q)=A_1(Q)^m$.
\end{proposition}

\begin{proof}
The first and second equalities follow from \cite[Theorem~1]{BereleEtAl2026}.
A nonzero product of $m$ matrix units from $A_1(Q)$ represents a
concatenation of $m$ paths of positive length. Hence
$A_1(Q)^m\subseteq A_m(Q)$. Conversely, dividing a path of length
at least $m$ into $m$ paths of positive length expresses its matrix
unit as a product of $m$ matrix units from $A_1(Q)$.
\end{proof}

The nonzero off-diagonal blocks of $A(Q)$ and $A_1(Q)$
are precisely the $(C,D)$-blocks with $C<_Q D$.
The corresponding block decompositions describe the products of
identity ideals used below.
In $A(Q)$, the diagonal block corresponding to an SCC $C$ is
$\Mat_{w(C)}(\kk)$.
In $A_1(Q)$, this block is unchanged when $C$ is cyclic and is a
$1\times1$ zero block otherwise.
In either algebra, the block from $C$ to $D$ is the full
rectangular matrix space when $C<_Q D$, and is zero when $C$ and
$D$ are incomparable. Indeed, a path between two SCCs can be
extended within them to join any given pair of their vertices.

\begin{definition}\label{def:chain-blocks}
Recall that $I_n=\Id(\Mat_n(\kk))$, and
let $I_0=\F$. For a nonempty tuple
$\mathbf a=(a_1,\cdots,a_r)$ of nonnegative integers, define
$
 I(\mathbf a)=I_{a_1}\cdots I_{a_r}$.
Let $d_i=\max\{a_i,1\}$.
The block algebra associated with $\mathbf a$ is the block upper
triangular subalgebra of $\Mat_{d_1+\cdots+d_r}(\kk)$ with
block sizes $d_1,\cdots,d_r$. Its $i$th diagonal block is
$\Mat_{a_i}(\kk)$ when $a_i>0$ and is zero when $a_i=0$.
Its $(i,j)$-block is $\Mat_{d_i\times d_j}(\kk)$ for $i<j$.
If all entries of $\mathbf a$ are positive, denote this algebra
by $U(\mathbf a)$.
\end{definition}

\begin{proposition}\label{prop:chain-block-identities}
For every nonempty tuple $\mathbf a$ of nonnegative integers,
the block algebra in Definition~\ref{def:chain-blocks} has
identity ideal $I(\mathbf a)$.
In particular, $\Id(U(\mathbf a))=I(\mathbf a)$ when all entries
are positive.
\end{proposition}

\begin{proof}
This follows from 
\cite[Theorem~1]{Lewin1974} or
\cite[Theorem~1.1]{Berele2012}.
\end{proof}

To compare products associated with positive integer
tuples, we use the following relation.

\begin{definition}\label{def:chain-order}
Let $\alpha=(a_1,\cdots,a_s)$ and
$\beta=(b_1,\cdots,b_t)$ be nonempty tuples of
positive integers. Write $\alpha\preceq\beta$ if
there is a nondecreasing map
$\phi:\{1,\cdots,s\}\to\{1,\cdots,t\}$ such that
$
\sum_{\phi(i)=j}a_i\leq b_j$
for $1\leq j\leq t$.
The map $\phi$ need not be surjective.
If $j$ is not in its image, the corresponding sum
is zero and the inequality is automatic.
\end{definition}

\begin{proposition}\label{prop:single-chain}
The relation $\preceq$ is transitive.
For nonempty tuples $\alpha,\beta$ of positive integers,
$
 I(\beta)\subseteq I(\alpha)
$ if and only if $
 \alpha\preceq\beta$.
If $\alpha\preceq\beta$, there is an algebra embedding
$U(\alpha)\hookrightarrow U(\beta)$, which need not preserve
the identity element.
\end{proposition}

\begin{proof}
Composing the nondecreasing maps in Definition~\ref{def:chain-order}
and summing their inequalities proves transitivity.
The equivalence is \cite[Theorem~2.3]{Berele2012}; it is also the
ordinary matrix case of Theorem~\ref{thm:mixed-capacity}.
If $\alpha=(a_1,\cdots,a_s)\preceq\beta$, choose a
nondecreasing map $\phi$ satisfying the inequalities
in Definition~\ref{def:chain-order}. Then choose, for every
$j$ and every $i$ with $\phi(i)=j$, disjoint sets of $a_i$ indices
in the $j$th diagonal block of $U(\beta)$.
Sending matrix units to the corresponding matrix units gives
the required embedding.
\end{proof}

We now associate two tuples of nonnegative integers with each
chain of SCCs to express the identity ideals of the path algebra
and its arrow ideal.

\begin{definition}\label{def:chain-weights}
A chain means a nonempty chain in $\mathcal P(Q)$ with respect
to $<_Q$. Let $\mathscr M(Q)$ be the set of maximal chains.
For $\mathcal C=(C_1<_Q\cdots<_Q C_r)$, define
\[
 \mathbf w(\mathcal C)=(w(C_1),\cdots,w(C_r)),
 \qquad
 a_i(\mathcal C)=
 \begin{cases}
 w(C_i),&C_i\text{ is cyclic},\\
 0,&C_i\text{ is not cyclic},
 \end{cases}
\]
and let $\mathbf a(\mathcal C)=(a_1(\mathcal C),\cdots,a_r(\mathcal C))$.
\end{definition}

The distinction between $\mathbf w(\mathcal C)$ and
$\mathbf a(\mathcal C)$ records exactly the diagonal matrix units
no longer present when trivial paths are removed.

\begin{proposition}\label{prop:chain-formulas}
Let $Q$ be a finite PI quiver. Then $
 \Id(\kk Q)
 =\bigcap_{\mathcal C\in\mathscr M(Q)}
     I(\mathbf w(\mathcal C))$ and
$ \Id(\kk Q_{\ge1})
 =\bigcap_{\mathcal C\in\mathscr M(Q)}
     I(\mathbf a(\mathcal C))$.
The same equalities also hold when the intersections are taken over all chains.
\end{proposition}

\begin{proof}
The formulas for intersections over maximal chains follow from Proposition~\ref{prop:tail-models}
and \cite[Theorem~1.1]{Berele2012}.
The equalities for all chains follow
from those for maximal chains.
\end{proof}

Several basic properties of arrow ideals follow directly.

\begin{proposition}\label{prop:acyclic-identities}
Suppose $Q$ is acyclic, and let $\ell$ be the maximum length of
a path in $Q$, allowing $\ell=0$. Then
$
 \Id(\kk Q_{\ge1})=\F^{\ell+1}
 =\langle x_1\cdots x_{\ell+1}\rangle_T$ and $
 \Id(\kk Q)=I_1^{\ell+1}$.
Moreover, $\Id(\kk Q)\subsetneq\Id(\kk Q_{\ge1})$.

In particular, for $n\ge2$, 
the arrow ideal $\kk(\overrightarrow{A_n})_{\ge1}$ of the linearly oriented quiver
$\overrightarrow{A_n}$ is isomorphic to the algebra of strictly
upper triangular $n\times n$ matrices, and both have identity ideal $\F^n$.
\end{proposition}

\begin{proof}
Every SCC is a noncyclic singleton. The maximum number of
vertices in a chain of the reachability order is $\ell+1$.
Proposition~\ref{prop:chain-formulas} gives the two ideal equalities;
the equality for $\kk Q$ is also proved in
\cite{CerulliIrelliEtAl2026}.
The ideal $\F^d$ is the span of all words of length at least $d$
and is generated as a $T$-ideal by $x_1\cdots x_d$.
The maximum length of a path in $\overrightarrow {A_n}$ is $n-1$, and $\kk(\overrightarrow {A_n})_{\ge1}$ is isomorphic to
the algebra of
strictly upper triangular $n\times n$ matrices, so they both have identity
ideal $\F^n$.
\end{proof}

A cyclic SCC gives a lower bound for the degree of a nonzero
identity of the arrow ideal.

\begin{proposition}\label{prop:cycle-degree}
If a finite PI quiver $Q$ has an SCC
that is a directed cycle with $n$ vertices,
then every nonzero polynomial identity of $\kk Q_{\ge1}$ has 
 degree at least $2n$.
\end{proposition}

\begin{proof}
The arrow ideal of the directed cycle is a subalgebra of
$\kk Q_{\ge1}$. Its matrix model is $\Mat_n(\kk)$. Hence every identity of $\kk Q_{\ge1}$ is an identity of
$\Mat_n(\kk)$.
The degree bound follows from Theorem~\ref{thm:amitsur-levitzki}.
\end{proof}

\subsection{Identities of path algebras and path ideals}
\label{sec:equality}

The inclusion $\kk Q_{\ge1}\subseteq\kk Q$ gives
$\Id(\kk Q)\subseteq\Id(\kk Q_{\ge1})$.
We first determine when the path algebra and its path ideals
have nonzero polynomial identities. 
For PI quivers, we use Proposition~\ref{prop:chain-formulas}
and Theorem~\ref{thm:finite-intersection} to obtain a criterion
for equality at $m=1$.
 We then prove that equality at any one
positive value of $m$ is equivalent to equality for every $m$.
Finally, we describe the stable matrix model and its identity ideal.

\begin{proposition}\label{prop:positive-pi}
Let $Q$ be a finite quiver. The following conditions are equivalent.
\begin{enumerate}
\item $\kk Q$ is a PI algebra.
\item $\kk Q_{\ge m}$ is a PI algebra for some $m\ge1$.
\item $\kk Q_{\ge m}$ is a PI algebra for every $m\ge1$.
\item Every SCC of $Q$ is either a directed
cycle or a single vertex without a loop.
\end{enumerate}
If these conditions fail, $\Id(\kk Q)=\Id(\kk Q_{\ge m})=0$
for every $m\ge1$.
\end{proposition}

\begin{proof}
The implications \textup{(i)}$\Rightarrow$\textup{(iii)}$\Rightarrow$\textup{(ii)}
follow from $\kk Q_{\ge m}\subseteq\kk Q$.
The equivalence of \textup{(i)} and \textup{(iv)} is shown in \cite{CerulliIrelliEtAl2026}.

Assume that $\kk Q_{\ge m}$  is a PI algebra for some $m \ge 1$.
If $\kk Q$ were not PI, 
then
\cite[Theorem~1.1]{CerulliIrelliEtAl2026} would give two 
closed paths $c_1$ and $c_2$ based at the same vertex such that $c_1c_2\ne c_2c_1$.
Choose two positive integers $i$ and $j$ such that $c_1^i$ and $c_2^j$ are of length at least $m$, and $c_1^ic_2^j \neq  c_2^j    c_1^i$.  
The nonunital subalgebra $B(c_1^i,c_2^j)$ of $\kk Q_{\ge m}$
generated by $c_1^i$ and $c_2^j$ is isomorphic to the free
nonunital algebra $\kk\langle x,y\rangle_+$.
Hence
$\Id(\kk Q_{\ge m})\subseteq\Id(B(c_1^i,c_2^j))=0$,
contrary to the assumption that $\kk Q_{\ge m}$ is PI.
Therefore $\kk Q$ is a PI algebra.

\end{proof}

For PI quivers, we now compare the identity ideals using
the relation $\preceq$ from Definition~\ref{def:chain-order}
and the criterion in Proposition~\ref{prop:single-chain}.
The next lemma uses the identity element to remove
the factors corresponding to zero entries.

\begin{lemma}\label{lem:delete-zero-factors}
Let $\mathbf a$ be a nonempty tuple of nonnegative integers, and let
$\mathbf a^+$ be the tuple obtained by deleting its zero entries.
Let $B$ be a nonzero unital algebra. Then $I(\mathbf a)\subseteq\Id(B)$ if and only if $\mathbf a^+$ is nonempty and
$I(\mathbf a^+)\subseteq\Id(B)$.
\end{lemma}

\begin{proof}
Suppose first that $I(\mathbf a)\subseteq\Id(B)$.
If every entry of $\mathbf a$ is zero, then $I(\mathbf a)$ is
a power of $\F$ and contains a monomial. Evaluating its variables
at $1_B$ contradicts $I(\mathbf a)\subseteq\Id(B)$.
Thus $\mathbf a^+$ is nonempty.

For each positive entry $a_i$, choose a polynomial $f_i\in I_{a_i}$,
using disjoint sets of variables for different entries. For each zero
entry, choose a new variable $z_i$. The product of these polynomials and variables lies in
$I(\mathbf a)$ and is an identity of $B$.
Evaluating every $z_i$ at $1_B$ gives
$I(\mathbf a^+)\subseteq\Id(B)$.

Conversely, $I(\mathbf a)\subseteq I(\mathbf a^+)\subseteq\Id(B)$.

\end{proof}

\begin{theorem}\label{thm:equality}

Let $Q$ be a finite PI quiver. Then
$\Id(\kk Q)=\Id(\kk Q_{\ge1})$ if and only if
for every maximal chain $\mathcal C$ in $\mathcal P(Q)$,
there is a chain $\mathcal D$ which is maximal among chains
consisting entirely of cyclic SCCs and satisfies
$\mathbf w(\mathcal C)\preceq\mathbf w(\mathcal D)$.

\end{theorem}

\begin{proof}
Suppose first that $\Id(\kk Q)=\Id(\kk Q_{\ge1})$, and let
$\mathcal C\in\mathscr M(Q)$. Proposition~\ref{prop:chain-formulas}
gives
$
\bigcap_{\mathcal L\in\mathscr M(Q)}
 I(\mathbf a(\mathcal L))
 \subseteq I(\mathbf w(\mathcal C))$.
Every entry of $\mathbf w(\mathcal C)$ is positive, so $I(\mathbf w(\mathcal C))$ is a product of proper verbally prime $T$-ideals.
Theorem~\ref{thm:finite-intersection} gives
$\mathcal L\in\mathscr M(Q)$ with
$
 I(\mathbf a(\mathcal L))
 \subseteq I(\mathbf w(\mathcal C))
 =\Id(U(\mathbf w(\mathcal C)))
$.
Since $U(\mathbf w(\mathcal C))$ is nonzero and unital,
Lemma~\ref{lem:delete-zero-factors} shows that deleting the
noncyclic SCCs of $\mathcal L$ gives a nonempty chain
$\mathcal D$ with
$
 I(\mathbf w(\mathcal D))
 \subseteq I(\mathbf w(\mathcal C))$.
Proposition~\ref{prop:single-chain} gives
$\mathbf w(\mathcal C)\preceq\mathbf w(\mathcal D)$.
Extend $\mathcal D$ to a chain maximal among chains consisting
entirely of cyclic SCCs. The same relation holds for the enlarged
chain, since the added SCCs may be left unused.

Conversely, suppose the chain condition holds.
For every $\mathcal C\in\mathscr M(Q)$, choose a chain
$\mathcal D$ as in the statement.
The matrix units of $A_1(Q)$ with both indices in the
SCCs of $\mathcal D$ span a subalgebra isomorphic to
$U(\mathbf w(\mathcal D))$.
By Proposition~\ref{prop:single-chain}, the condition
$\mathbf w(\mathcal C)\preceq\mathbf w(\mathcal D)$ gives an embedding
$U(\mathbf w(\mathcal C))\hookrightarrow U(\mathbf w(\mathcal D))$. Hence
$ \Id(\kk Q_{\ge1})
 \subseteq I(\mathbf w(\mathcal D))
 \subseteq I(\mathbf w(\mathcal C))$.
Intersecting over $\mathcal C\in\mathscr M(Q)$ and applying
Proposition~\ref{prop:chain-formulas} gives
$\Id(\kk Q_{\ge1})\subseteq\Id(\kk Q)$.
The reverse inclusion follows from $\kk Q_{\ge1}\subseteq\kk Q$.
\end{proof}

\begin{remark}
If every SCC of $Q$ is cyclic, then
$\Id(\kk Q)=\Id(\kk Q_{\ge1})$, since one may take
$\mathcal D=\mathcal C$ in Theorem~\ref{thm:equality}.
In general, the chains $\mathcal C$ and $\mathcal D$
need not lie in the same connected component of the
underlying undirected graph.
\end{remark}

A path through a cyclic SCC can be lengthened without changing
its starting or terminal vertex.
Consequently, for every cyclic chain $\mathcal D$,
each $A_m(Q)$ contains a subalgebra isomorphic to
$U(\mathbf w(\mathcal D))$.

\begin{corollary}\label{cor:all-tail-equality}
Let $Q$ be a finite quiver. If $Q$ is PI, the following
conditions are equivalent. 
\begin{enumerate}
\item $\Id(\kk Q)=\Id(\kk Q_{\ge m})$ for some $m\ge1$.
\item $\Id(\kk Q)=\Id(\kk Q_{\ge m})$ for every $m\ge1$.
\item For every maximal chain $\mathcal C$ in $\mathcal P(Q)$,
there is a chain $\mathcal D$ which is maximal among chains
consisting entirely of cyclic SCCs and satisfies
$\mathbf w(\mathcal C)\preceq\mathbf w(\mathcal D)$.
\end{enumerate}
If $Q$ is not PI, $\Id(\kk Q)=\Id(\kk Q_{\ge m})=0$ for every $m\ge1$.

\end{corollary}

\begin{proof}
Suppose first that $Q$ is PI.
If equality holds for some $m$, then
$
 \Id(\kk Q)\subseteq\Id(\kk Q_{\ge1})
 \subseteq\Id(\kk Q_{\ge m})=\Id(\kk Q)$.
Theorem~\ref{thm:equality} therefore gives \textup{(iii)}. \textup{(ii)} implies \textup{(i)}.

Now assume that \textup{(iii)} holds. Let $m\ge1$.
For each $\mathcal C\in\mathscr M(Q)$, choose $\mathcal D$
as in \textup{(iii)}.
 $A_m(Q)$ contains a subalgebra isomorphic to
$U(\mathbf w(\mathcal D))$.
Indeed, each corresponding matrix unit represents a path
which meets a cyclic SCC, and repeating a closed path gives
such a path of length at least $m$.  Hence
$
 \Id(\kk Q_{\ge m})
 \subseteq I(\mathbf w(\mathcal D))
 \subseteq I(\mathbf w(\mathcal C))$.
Intersecting over $\mathcal C\in\mathscr M(Q)$ and applying
Proposition~\ref{prop:chain-formulas} gives
$\Id(\kk Q_{\ge m})\subseteq\Id(\kk Q)$.
The reverse inclusion follows from
$\kk Q_{\ge m}\subseteq\kk Q$.
The non-PI case follows from Proposition~\ref{prop:positive-pi}.
\end{proof}

The corollary concerns equality of the identities of the path
algebra and its path ideals. For a finite PI quiver, the descending
sequence of matrix models $A_m(Q)$ stabilizes even when this
equality fails. We now describe its stable term and bound
the value of $m$ at which it is reached.

\begin{proposition}\label{prop:tail-stabilization}
Let $Q$ be a finite PI quiver. Let $\ell_{\mathrm{nc}}$ be the maximum
length of a path all of whose vertices lie in noncyclic SCCs,
including trivial paths. If there are no such vertices, let
$\ell_{\mathrm{nc}}=-1$. Define
$H(Q)=\max\{1,\ell_{\mathrm{nc}}+1\}$, and let
$e=\sum_{v}E_{vv}$, where the sum ranges over the vertices in cyclic
SCCs. Let $A_\infty(Q)=A_1(Q)eA_1(Q)$. For every $m\ge H(Q)$,
$
 A_m(Q)=A_\infty(Q)$ and $
 \Id(\kk Q_{\ge m})=\Id(A_\infty(Q))$.
The algebra $A_\infty(Q)$ is spanned by the matrix units $E_{uv}$
for which a path from $u$ to $v$ meets a cyclic SCC.
If $n_{\mathrm{nc}}$ is the number of vertices in noncyclic SCCs,
then $H(Q)\le\max\{1,n_{\mathrm{nc}}\}$.
\end{proposition}

\begin{proof}

Expanding $e$ gives
$
A_\infty(Q)=\operatorname{span}_{\kk}
\{E_{uv}\mid E_{uc},E_{cv}\in A_1(Q)
\text{ for some }c\text{ in a cyclic SCC}\}.
$
Since $E_{cc}\in A_1(Q)$ for every vertex $c$ in a cyclic SCC,
the condition in this expression is equivalent to the existence
of a path from $u$ to $v$ meeting a cyclic SCC.
This proves the asserted spanning description.

Every path meeting a cyclic SCC can be lengthened arbitrarily
by repeating a closed path, without changing its endpoints.
Hence every matrix unit in this spanning set belongs to
$A_m(Q)$ for every $m\ge1$, so $A_\infty(Q)\subseteq A_m(Q)$.
If $m\ge H(Q)$, then $m>\ell_{\mathrm{nc}}$, and every path
of length at least $m$ must meet a cyclic SCC.
The spanning description therefore gives
$A_m(Q)\subseteq A_\infty(Q)$.
Consequently, $A_m(Q)=A_\infty(Q)$ for $m\ge H(Q)$.
By Proposition~\ref{prop:tail-models},
$\Id(\kk Q_{\ge m})=\Id(A_m(Q))=\Id(A_\infty(Q))$
for $m\ge H(Q)$.

If $n_{\mathrm{nc}}>0$, a path whose vertices lie in
noncyclic SCCs cannot repeat a vertex, so its length is
at most $n_{\mathrm{nc}}-1$.
Hence $\ell_{\mathrm{nc}}\le n_{\mathrm{nc}}-1$ and
$H(Q)\le n_{\mathrm{nc}}$.
If $n_{\mathrm{nc}}=0$, then $H(Q)=1$ by definition.
\end{proof}

\begin{remark}\label{rem:stable-nonunital-blocks}
The bound in Proposition~\ref{prop:tail-stabilization} is attained
by the linearly oriented quiver
$\overrightarrow{A_n}$. Indeed, $A_m(\overrightarrow{A_n})$
is nonzero for $1\le m<n$ and is zero for $m\ge n$.

The algebra $A_\infty(Q)$ may contain $E_{uv}$ even when
$u$ or $v$ lies in a noncyclic SCC.
For example, let $Q$ consist of a loop at vertex $1$
and an arrow from $1$ to vertex $2$. Then
\[
 A_\infty(Q)=A_1(Q)=
 \left\{
 \begin{pmatrix}a&b\\0&0\end{pmatrix}
 \middle| a,b\in\kk
 \right\}.
\]
Here $e=E_{11}$, so $eA_1(Q)e=\kk E_{11}\cong\kk$.
The identity ideal of $A_\infty(Q)$ is $I_1\F$,
whereas the identity ideal of $eA_1(Q)e$ is $I_1$.
In particular, $[E_{11},E_{12}]=E_{12}\ne0$
in $A_\infty(Q)$.
\end{remark}

\subsection{Numerical criteria and examples}
\label{sec:examples}

Throughout this subsection, let $Q$ be a finite PI quiver.
When all cyclic SCCs have the same number of vertices, we
express the equality criterion in Theorem~\ref{thm:equality}
as a numerical condition on the SCCs in each chain.
We illustrate the criterion with examples, including one in
which an explicit polynomial is an identity of the arrow ideal
but not of the path algebra.

\begin{proposition}\label{prop:uniform-cycles}
Let $Q$ be a finite PI quiver with at least one cyclic SCC,
and suppose every cyclic SCC has $m$ vertices.
Let $c$ be the maximum number of SCCs in a chain consisting
entirely of cyclic SCCs.
For each maximal chain $\mathcal C\in\mathscr M(Q)$ containing
$k$ cyclic SCCs, let $t_0,\cdots,t_k$ be the numbers of noncyclic
SCCs before, between, and after these cyclic SCCs, in their order
along $\mathcal C$.
When $k=0$, let $t_0$ be the number of SCCs in $\mathcal C$.
Define
$\nu_m(\mathcal C)
=k+\sum_{i=0}^{k}\lceil t_i/m\rceil$.
Then $\Id(\kk Q)=\Id(\kk Q_{\ge1})$ if and only if
$\max_{\mathcal C\in\mathscr M(Q)}\nu_m(\mathcal C)=c$.
\end{proposition}

\begin{proof}
Let $\mathcal C\in\mathscr M(Q)$, and let $\mathcal D$
be a chain consisting of $s$ cyclic SCCs.
Then $\mathbf w(\mathcal D)=(m,\cdots,m)$.

To satisfy the inequalities in Definition~\ref{def:chain-order},
each cyclic SCC of $\mathcal C$ must be assigned to a distinct
SCC of $\mathcal D$, and no other SCC of $\mathcal C$ can be
assigned to the same one.
Each of the $t_i$ consecutive noncyclic SCCs $C$ satisfies
$w(C)=1$, so they require at least $\lceil t_i/m\rceil$
additional SCCs of $\mathcal D$.
Since the assignment is nondecreasing, these additional SCCs
are distinct for different values of $i$.
These numbers suffice, since one can assign at most $m$
successive noncyclic SCCs to each additional SCC of $\mathcal D$.
Thus
$
 \mathbf w(\mathcal C)\preceq\mathbf w(\mathcal D)
$ if and only if $
 \nu_m(\mathcal C)\le s$.

By the definition of $c$ and Theorem~\ref{thm:equality},
$\Id(\kk Q)=\Id(\kk Q_{\ge1})$ if and only if
$\nu_m(\mathcal C)\le c$ for every
$\mathcal C\in\mathscr M(Q)$.
A chain consisting of $c$ cyclic SCCs extends to a maximal
chain $\mathcal C$, for which $\nu_m(\mathcal C)\ge c$.
Hence
$\max_{\mathcal C\in\mathscr M(Q)}\nu_m(\mathcal C)\ge c$,
and the asserted equivalence follows.
\end{proof}

\begin{corollary}\label{cor:loops-height}
Suppose every cyclic SCC of $Q$ is a loop.
Let $h$ be the maximum number of SCCs in a chain, and let $c$ be
the maximum number of cyclic SCCs in a chain.
Then $\Id(\kk Q)=\Id(\kk Q_{\ge1})$ if and only if $h=c$.
Equivalently, equality holds if and only if there is a longest
chain consisting entirely of cyclic SCCs.
\end{corollary}

\begin{proof}
If $c>0$, apply Proposition~\ref{prop:uniform-cycles} with $m=1$,
since $\nu_1(\mathcal C)$ is the number of SCCs in $\mathcal C$.
If $c=0$, then $h>0$, and the two identity ideals are distinct
by Proposition~\ref{prop:acyclic-identities}.
\end{proof}

The equality may also hold when $Q$ has noncyclic SCCs,
as the following example shows.

\begin{example}\label{ex:equal-size-capacity}
\[
Q_+\;
\left\{
\begin{gathered}
\begin{tikzcd}[
  ampersand replacement=\&,
  column sep=1.5em,
  row sep=1.5em
]
a_1 \arrow[r,bend left=20,"\alpha_1"]
\& a_2 \arrow[l,bend left=20,"\alpha_2"]
        \arrow[r,"\alpha_3"]
\& a_3 \arrow[r,bend left=20,"\alpha_4"]
\& a_4 \arrow[l,bend left=20,"\alpha_5"]
\\[1em]
a_5 \arrow[r,"\alpha_6"]
\& a_6 \arrow[r,"\alpha_7"]
\& a_7 \arrow[d,bend left=20,"\alpha_8"]
\&
\\
\& \& a_8 \arrow[u,bend left=20,"\alpha_9"] \&
\end{tikzcd}
\end{gathered}
\right.
\qquad
\left.
\begin{gathered}
\begin{tikzcd}[
  ampersand replacement=\&,
  column sep=1.5em,
  row sep=1.5em
]
b_1 \arrow[r,bend left=20,"\beta_1"]
\& b_2 \arrow[l,bend left=20,"\beta_2"]
        \arrow[r,"\beta_3"]
\& b_3 \arrow[r,bend left=20,"\beta_4"]
\& b_4 \arrow[l,bend left=20,"\beta_5"]
\\[1em]
b_5 \arrow[r,"\beta_6"]
\& b_6 \arrow[r,"\beta_7"]
        \arrow[d,bend left=20,"\beta_8"]
\& b_7
\&
\\
\& b_8 \arrow[u,bend left=20,"\beta_9"] \& \&
\end{tikzcd}
\end{gathered}
\right\}
\;Q_-
\]
Let $Q_+$ and $Q_-$ be the quivers on the left and right,
respectively, in the figure above.
In both quivers, a longest chain consisting entirely
of cyclic SCCs has two SCCs.

The two maximal chains $\mathcal C_1^+,\mathcal C_2^+$
in $\mathcal P(Q_+)$ satisfy
$\mathbf w(\mathcal C_1^+)=(2,2)$ and
$\mathbf w(\mathcal C_2^+)=(1,1,2)$.
Since $\nu_2(\mathcal C_1^+)=\nu_2(\mathcal C_2^+)=2$,
Proposition~\ref{prop:uniform-cycles} gives
$\Id(\kk Q_+)=\Id(\kk(Q_+)_{\ge1})$.
The two maximal chains $\mathcal C_1^-,\mathcal C_2^-$
in $\mathcal P(Q_-)$ satisfy
$\mathbf w(\mathcal C_1^-)=(2,2)$ and
$\mathbf w(\mathcal C_2^-)=(1,2,1)$.
Here $\nu_2(\mathcal C_1^-)=2$ and
$\nu_2(\mathcal C_2^-)=3$.
Thus $\Id(\kk Q_-)\subsetneq\Id(\kk(Q_-)_{\ge1})$.
For $Q_+$, $(1,1,2)\preceq(2,2)$,
but for $Q_-$, $(1,2,1)\npreceq(2,2)$.
Indeed, in $(1,2,1)$ the entry $2$ must have
a different image from both entries $1$, and a nondecreasing
map cannot send these two entries $1$ to the same value.
\end{example}

The next example gives an explicit polynomial identity of
the arrow ideal which is not an identity of the path algebra.

\begin{example}\label{ex:cycle-tail}
\[
\begin{tikzcd}[column sep=large]
v_1 \arrow[r,"\alpha_1"]
& v_2 \arrow[r,"\alpha_2"]
& \cdots \arrow[r,"\alpha_{n-1}"]
& v_n \arrow[lll,bend left=35,"\alpha_n"]
      \arrow[r,"\beta"]
& v
\end{tikzcd}
\]
Let $Q$ be the quiver in the figure above.
For $n=1$, the directed cycle is a loop.
By Proposition~\ref{prop:chain-formulas},
$\Id(\kk Q)=I_nI_1$ and $\Id(\kk Q_{\ge1})=I_n\F$.
Let
$
 \st_{2n}$
be the standard polynomial of degree $2n$. Then
$
 \st_{2n}(x_1,\cdots,x_{2n})y
 \in\Id(\kk Q_{\ge1})\setminus\Id(\kk Q)$.
\end{example}

\begin{proof}
Theorem~\ref{thm:amitsur-levitzki} gives
$\st_{2n}\in I_n$, so
$\st_{2n}y\in I_n\F$.
In $U(n,1)$, evaluate $x_{2j-1}$ at $E_{jj}$ and
$x_{2j}$ at $E_{j,j+1}$ for $1\le j\le n$,
and evaluate $y$ at $E_{n+1,n+1}$.
Only the monomial $x_1\cdots x_{2n}$ in $\st_{2n}$
contributes a nonzero value, so the value of $\st_{2n}y$
is $E_{1,n+1}$.
Since $\Id(U(n,1))=I_nI_1=\Id(\kk Q)$, this polynomial
is not an identity of $\kk Q$.
\end{proof}

\bibliographystyle{plain}
\bibliography{Tideal}

\end{document}